\documentclass[preprint, authoryear,12pt]{elsarticle}

\usepackage[utf8]{inputenc}

\usepackage{amsmath,amssymb,amsthm,amscd}
\usepackage{mathtools}
\usepackage{dsfont} 
\usepackage{bm} 
\usepackage[binary-units=true]{siunitx}
\usepackage{multirow}
\usepackage{caption}
\usepackage[font=footnotesize]{subcaption}
\usepackage{booktabs} 
\usepackage{blkarray} 
\usepackage{pdflscape}
\usepackage{longtable}

\usepackage{tikz}
\usepackage{xcolor}
\usepackage{graphicx} 
\usepackage[colorlinks=true, citecolor=blue, linkcolor=blue]{hyperref}

\usepackage{paralist} 
\usepackage{verbatim} 

\usepackage{fancybox}
\usepackage{varioref} 
\usepackage{chapterbib} 
\usepackage{eurosym}
\usepackage[running]{lineno}
\usepackage{afterpage}

\usepackage[compatible]{nomencl}
\usepackage{xspace}
\newcommand{\dd}{\mathrm{d}} 

\newcommand{\norm}[1]{\left\lVert#1\right\rVert}

\newcommand{\pbk}[1]{\left(#1\right)}
\newcommand{\cbk}[1]{\left\lbrace#1\right\rbrace}
\newcommand{\sbk}[1]{\left\lbrack#1\right\rbrack}
\newcommand{\ie}{\textit{i.e.}\@\xspace}
\newcommand{\eg}{\textit{e.g.}\@\xspace}
\newcommand{\HE}{\hbox{H\kern-.12em\lower.48ex\hbox{E}}}
\def\leaderfill{\leaders\hbox to 1em{\hss.\hss}\hfill}
\DeclareSIUnit\molar{\mole\per\cubic\metre}
\DeclareSIUnit\Molar{\textsc{M}}

\graphicspath{ {images/} }

\journal{Journal of Chromatography A}
\usepackage{lineno}

\begin{document}

\begin{frontmatter}

\title{Clustering-based convergence diagnostic for multi-modal identification in parameter estimation of chromatography model with parallel MCMC}

\author[add1]{Yue-Chao Zhu}
\author[add2]{Zhaoxi Sun}
\author[add3]{Qiao-Le He\corref{cor1}}
\ead{qiaole.he@rwth-aachen.de}
\cortext[cor1]{Corresponding author.}

\address[add1]{State Key Laboratory of Bioreactor Engineering, East China University of Science and Technology, 200237 Shanghai, China}
\address[add2]{Beijing National Laboratory for Molecular sciences, Institute of Theoretical and Computational Chemistry, College of Chemistry and Molecular Engineering, Peking University, 100871 Beijing, China}
\address[add3]{Enzymaster (Ningbo) Bio-Engineering Co., Ltd., North Century Avenue 333, 315100 Ningbo, China}

\begin{abstract}
Uncertainties from experiments and models render multi-modal difficulties in model calibrations.
Bayesian inference and \textsc{mcmc} algorithm have been applied to obtain posterior distributions of model parameters upon uncertainty.
However, multi-modality leads to difficulty in convergence criterion of parallel \textsc{mcmc} sampling chains.
The commonly applied $\widehat{R}$ diagnostic does not behave well when multiple sampling chains are evolving to different modes.
Both partitional and hierarchical clustering methods has been combined to the traditional $\widehat{R}$ diagnostic to deal with sampling of target distributions that are rough and multi-modal.
It is observed that the distributions of binding parameters and pore diffusion of particle parameters are multi-modal.
Therefore, the steric mass-action model used to describe ion-exchange effects of the model protein, lysozyme, on the \textsc{sp} Sepharose \textsc{ff} stationary phase might not be fully capable in certain experimental conditions, as model uncertainty from steric mass-action would result in multi-modality.

\end{abstract}


\end{frontmatter}


\section{Introduction}
Model-based methods can be an alternative to experiments in mechanism exploit, design of experiments, process design and prediction \citep{guiochon2003modeling}.
Chromatography is one of the most commonly applied technologies in the downstream processing of biomoecules \citep{carta2010protein}.
In chromatographic field, model-based studies have emerged since 1940s \citep{martin1941new, james1952gas, wilson1940present, devault1943theory, craig1944identification} and, hitherto, extensive papers have been published \citep{guiochon2006fundamentals, antos2001application, aumann2007continuous, hahn2014adjoint, leweke2016fast}.
Identification of model parameters in model calibration step (\ie, parameter estimation) is inevitable in the framework of model-based studies, to further exploit models.
Deterministic \citep{vonLieres2010fast, hahn2014adjoint, forssen2006improved, osberghaus2012determination}, heuristic \citep{gao2005estimation, kaczmarski2006use, xu2013determination, huuk2014model}, and stochastic \citep{he2019modelbased, briskot2019prediction} algorithms have been adopted to solve the inverse problems of chromatographic models.

Recently, Bayesian inference has been applied to deal with the parameter estimation of chromatographic models, which generally renders versatile information \citep{briskot2019prediction, he2019bayesian}.
Particularly, Bayesian inference involves an interpretation of uncertainties, ranging from experiment, modelling, and numerical solution, into parameter description of probabilities.
Measurement of the uncertainties are propagated in a mathematically consistent manner.
The Bayesian perspective differs fundamentally with the frequentist one (\eg, maximum likelihood estimation), in which point estimate is attained. 

Markov Chain Monte Carlo (\textsc{mcmc}) is a key computational tool in Bayesian statistics.
Computationally, \textsc{mcmc} makes tactable many applied Bayesian analysis that were analytically intactable and impractical, by generating a sample of posterior draws.
It has proven to converge to the target distribution when the number of draws approaches infinity.
Although the applicability of \textsc{mcmc}, approaching \emph{infinity} could raise computational cost even in the computing world with exponential growth in processor speeds, when solving complex and higher-dimensional problems, in particular, problems that are full of many local optima.
In other words, two of the most common causes of slow convergence are multi-modality and dimensional dependence \citep{VanDerwerken2015monitoring}.
In the chromatography field, we are often fitting models with large numbers of parameters. 
Moreover, the strong non-linearity of the likelihood in the chromatographic field could cause multi-modality.
However, according to our limited knowledge, few studies have reported in the published literature and elucidated the multi-modality situations in parameter estimation of chromatographic models.
In practical implementation with insufficient sampling length, according to experience, it could render arbitrary bad results in some scenarios \cite{}.
Therefore, determining how \emph{finite} long to run is sufficient to mimic and reproduce the target distribution is an important and challenging topic.
Practically, we rely heavily on statistical diagnostics to judge and claim convergence. 
Various convergence diagnostics are available, \eg, \cite{heidelberger1983simulation}, \cite{cowles1996markov}, \cite{brooks1998general}, \cite{mengersen1999mcmc}, \cite{gelman1992inference}, \cite{raftery1992practical}, \cite{plummer2006coda}, \cite{paul2012assessing}.

Various technologies have been presented to enhance the convergence rate \citep{haario2006dram, braak2006markov}.
Parallelism is a strategy to speed-up convergence and identify multiple modes.
Due to its iterative nature, \textsc{mcmc} is inherently serial and thus can not be parallelized in programming (\eg, openMP). 
Parallelism can be applied in such ways: 1) multiple chains of \textsc{mcmc} are run in parallel, where the chains can be either independent (\ie, embarrassingly parallel \textsc{mcmc}) or communicative (interacting parallel \textsc{mcmc}).
Running multiple chains in parallel does not lead to significant speedups when one chain's convergence is slow, since each of them must be run long enough to equilibrate \citep{geyer1992practical}; this is also related to convergence diagnostic for multiple chains.
So the above parallelism idea seems contradictory in the literature (cf.~``many short runs'' versus ``one long run'' in \cite{geyer1991markov} for detailed information).
Since the multi-modal feature of target distribution in chromatography, ``many short runs'' are inclined in this study. 

Furthermore, running multiple chains is critical to any convergence diagnostics.
Among the various convergence diagnostics (cf.~\cite{cowles1996markov}, \cite{mengersen1999mcmc}, \cite{robert2011short} for review), the potential scale reduction factor, $\widehat{R}$, is the most widely used one.
It has been integrated into widely applied software packages (\eg, Stan, JAGS, WinBUGS, openBUGS, PyMC3).
The factor $\widehat{R}$ starts several Markov chains at overdispersed initial values, and monitors convergence by comparing between- and within-chain variances for selected scalar functions of the chains by means of a scale reduction factor. 
A multivariate extension introduced later by \cite{brooks1998general}. But, $\widehat{R}$ still fails to diagnose poor mixing when the chain has a heavy tail or when the variance varies across the chains. 
\cite{vehtari2021rank} then introduced rank-based split-$\widehat{R}$ to fix these problems.

However, the $\widehat{R}$ diagnostic is still not without limitations. 
In sampling of rough or multi-modal posterior distribution in chromatographic field, where the above mentioned parallelism technique has been applied to boost, the $\widehat{R}$ diagnostic does not behave well.
Specifically, in order to have $\widehat{R}$, \eg, $< 1.10$, each chain has to explore all the modes such that variations between chains are small and equilibrated; but Markov chains can get trapped in a local mode with extremely low possibilities for transitions between modes.
In these scenarios, $\widehat{R}$ could keep unchanged for an extremely long time.

In this study, we present $K$-means clustering based $\widehat{R}$ diagnostic to assess convergence of respectively partitioned groups (\ie, each group for each mode) of multiple chains of \textsc{mcmc}.
The posterior distribution is partitioned into $K$ sub-posterior distributions, $p(\theta|y) = \prod_{i=1}^K \omega_i p_i(\theta|y)$.
For each sub-posterior distribution $p_i$, which is assumed to be unimodal, \textsc{mcmc} sampling can be boosted, as random-walk is able to converge quickly to a local mode of the target distribution; and the traditional $\widehat{R}$ is able to diagnose convergence.
Both partitional and hierarchical clustering methods shall be implemented anc compared.
Due to its simplicity, robustness and speed, $K$-means is one of the most popular and widely used methods of clustering.
$K$-means clustering is an expression of alternate minimization that aims at partitioning $n$ observations (here, multiple chains) into $K$ clusters, in which each observation belongs to the cluster with the nearest mean \citep{lloyd1982least}.
It utilizes a centroid-based approach for clustering. 
In specific, we specify the number of clusters ($K$) and the method then identifies which chains belong to which cluster.
$K$-means++ initialization is further adopted to alleviate a drawback of the $K$-means method --- randomly generated centroids can be arbitrary bad.
Other variant methods, such as $K$-median, $K$-medoids, are applied to validate the clustering of multiple chains.
In most cases, the $K$ is determined in a heuristic way. 
Hierarchical clustering, such as, in dendrogram representation could be another option without finding the optimal $K$ value.
After clustering of multiple chains, we do not need to wait extreme long time to let each chain explore all the modes, such that $\widehat{R}$ diagnostic behaves well in each cluster.
Further, delayed rejection and adaptive Metropolis strategies are used to enhance efficiency of the na\"ive \textsc{mcmc} sampling.
Capability of \textsc{mcmc} to identify multiple modes, along with convergence diagnostic to stop sampling, Bayesian inference can be better utilized in the chromatography field to render versatile information, such as, uncertainty of binding modes on parameters values that have not been reported.

\section{Column models}
\begin{figure}[h]
    \centering
    \includegraphics[width=0.8\textwidth]{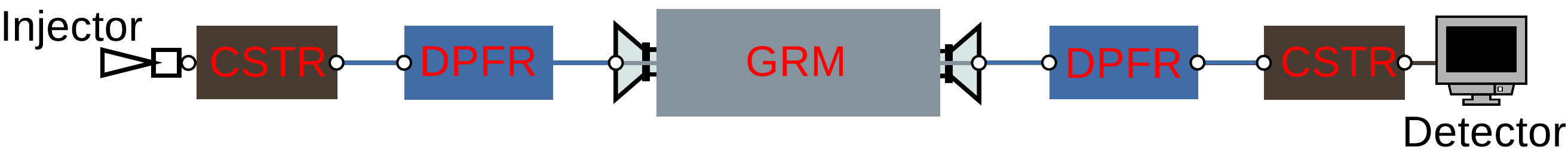}
    \caption{Schematic of mathematical models used to describe transport behaviour between injector and detector of a chromatography equipment. Lines do not represent tubings, instead just connection. External dispersion effects between injector and column are described with \textsc{dpfr}; external volume effects are described with \textsc{cstr}. Mass transfer in the column is described with \textsc{grm}. Models between column and detector are mirrored in this study.}
    \label{fig:model_scheme}
\end{figure}

A serial of mathematical models are concatenated in this study to describe the whole transport behaviour of studied components between the injector and the detector of a chromatography equipment, as shown in Fig.~\ref{fig:model_scheme}. 
External dispersion effects of the system volume between the injector and the column are described with a dispersive plug flow reactor (\textsc{dpfr}); external volume shifts are described with a continuous stirred tank reactor (\textsc{cstr}). 
Lines in Fig.~\ref{fig:model_scheme} do not represent physical tubings, instead just connector of models. 
Various levels of mass transfer resistance in the chromatogrpahic column are described with general rate model (\textsc{grm}).
Models between the column and the detector are mirrored.

\subsection{CSTR and DPFR}
The \textsc{cstr} is described by:
\begin{equation}
    \frac{\partial \hat{c}_i}{\partial t} = \frac{Q}{V^\textsc{cstr}} \pbk{ \hat{c}^\text{in}_i - \hat{c}_i } 
    \label{eq:cstr}
\end{equation}
And the \textsc{dpfr} involves convection and dispersion contributions:

\begin{equation}
    \frac{\partial \tilde{c}_i}{\partial t} = - u^\textsc{dpfr} \frac{\partial \tilde{c}_i}{\partial z} + D^\textsc{dpfr}_\text{ax} \frac{\partial^2 \tilde{c}_i}{\partial z^2} 
    \label{eq:dpfr}
\end{equation}
with inlet and outlet boundary conditions:

\begin{subequations}
    \begin{align}
        \sbk{ u^\textsc{dpfr} \tilde{c}_i - D^\textsc{dpfr}_\text{ax} \frac{\partial \tilde{c}_i}{\partial z} }_{z=0} & = u^\textsc{dpfr}\, \tilde{c}^\text{in}_i\\
                \left. \frac{\partial \tilde{c}_i}{\partial z} \right|_{z = L^\textsc{dpfr}} & = 0 
    \end{align}
\end{subequations}
where $\tilde{c}_i$ and $\hat{c}_i$ are the component concentrations of units \textsc{dpfr} and \textsc{cstr}.
$\hat{c}_i^\text{in}$ is the injection profile; the inlet profile of the \textsc{dpfr} unit is the outlet profile of the previous \textsc{cstr} unit, $\tilde{c}_i^\text{in} = \hat{c}_i^\text{out}$.
Moreover, $u^\textsc{dpfr}$, $D^\textsc{dpfr}_\text{ax}$ and $L^\textsc{dpfr}$ denote linear velocity, axial dispersion and length of the \textsc{dpfr}, while $Q$ and $V^\textsc{cstr}$ are the volumetric flow rate and the volume of the \textsc{cstr}.

\subsection{General rate model}
The \textsc{grm} is used to describe the transport behaviour of components in the chromatographic column.
Convection and axial dispersion in the bulk liquid are considered, as well as film mass transfer, pore diffusion and surface diffusion in the porous beads:

\begin{subequations} \label{eq:GRM}
    \begin{align}
    \frac{\partial c_i}{\partial t} & = -u_\text{int} \frac{\partial c_i}{\partial z} + D_{\text{ax}} \frac{\partial^2 c_i}{\partial z^2} - \frac{1-\varepsilon_c}{\varepsilon_c} \frac{3}{r_p} k_{f,i} \left(c_i - c_{p,i}(r\!=\!r_p)\right) \\
    \frac{\partial c_{p,i}}{\partial t} & = D_{p,i} \frac{1}{r^2} \frac{\partial}{\partial r} \left(r^2 \frac{\partial c_{p,i}}{\partial r} \right) + D_{s,i} \frac{1}{r^2} \frac{\partial}{\partial r} \left(r^2 \frac{\partial q_i}{\partial r} \right) - \frac{1-\varepsilon_p}{\varepsilon_p}\frac{\partial q_i}{\partial t}
    \end{align}
\end{subequations}
In Eq.~\eqref{eq:GRM}, $z \in [0,L]$ denotes the axial position where $L$ is the column length, while $r \in [0,r_p]$ denotes the radial position where $r_p$ is the particle radius.
Furthermore, $c_i$, $c_{p,i}$ and $q_i$ denote the interstitial, stagnant and stationary phase concentrations of component $i\in \{1,\dots,m\}$ in the column, respectively.
$t$ is the time, $\varepsilon_c$ and $\varepsilon_p$ are the column and particle porosities, $u_{\text{int}}$ interstitial velocity; $D_{\text{ax}}$ is the axial dispersion coefficient, $D_{p,i}$ the effective pore diffusion coefficient, $D_{s,i}$ the surface diffusion coefficient, and $k_{f,i}$ the film mass transfer coefficient.
At the column inlet ($z=0$) and outlet ($z=L$), Danckwerts boundary conditions \citep{Barber1998Boundary} are applied:
\begin{equation}
    \begin{split}
        \sbk{ u_{\text{int}} \, c_i - D_{\text{ax}} \dfrac{\partial c_i}{\partial z} }_{z=0} &= u_{\text{int}} \, c^\text{in}_i \\
        \left. \dfrac{\partial c_i}{\partial z} \right|_{z=L} &= 0
    \end{split}
\end{equation}
where $c^\text{in}_i = \tilde{c}_i^\text{out}$ is the inlet concentration of component $i$.
The boundary conditions at the particle surface ($r = r_p$) and centre ($r = 0$) are described by:

\begin{equation}
    \begin{split}
        \sbk{ k_{f,i}\, c_{p,i} - \varepsilon_p D_{p,i} \dfrac{\partial c_{p,i}}{\partial r} }_{r=r_p} &= k_{f,i}\, c_{p,i} (r\!=\!r_p) \\
        \left. \dfrac{\partial c_{p,i}}{\partial r} \right|_{r=0}  &= 0
    \end{split}
    \label{eq:Danckwerts_surface}
\end{equation}

\subsection{Steric mass-action isotherm}
The change rate of the stationary phase concentration, $\frac{\partial q_i}{\partial t}$, in Eq.~\eqref{eq:GRM} is calculated with binding isotherm models.
Ion-exchage effects of protein molecules to the stationary phase can be well-accounted with the steric mass-action (\textsc{sma}) model \citep{brooks1992steric}.
The thermodynamic binding can be considered either in a dynamic form, see Eq.~\eqref{eq:SMA_dyn}, or an implicit algebraic form when the binding is quasi-stationary, see Eq.~\eqref{eq:SMA_alg}.

\begin{subequations} \label{eq:SMA}
    \begin{align} 
        \frac{\partial q_i}{\partial t} &= k_{a,i}\, c_{p,i}\, \bar{q}_0^{\nu_i} - k_{d,i}\, q_i\, c_{p,0}^{\nu_i} \label{eq:SMA_dyn} \\
        q_i & = k_{\text{eq},i}\, c_{p,i} \pbk{ \frac{ \bar{q}_0 }{c_{p,0}} }^{\nu_i} \label{eq:SMA_alg}
    \end{align}
\end{subequations}
where $k_{\text{eq},i} = \frac{k_{a,i}}{k_{d,i}}$ is the equilibrium constant, $\nu$ denotes the characteristic charge, $\sigma$ the shielding factor, and $\Lambda$ the ionic capacity.
Component $i=0$ is introduced to denote the salt component, such that $c_{p,0}$ and $\bar{q}_0 =  \Lambda - \sum_{j=1}^m (\nu_j + \sigma_j)\, q_j$ represent salt concentration in the porous stagnant phase and salt concentration that available for binding in the stationary phase.

\subsection{Initial conditions}
Apart from boundary conditions, initial conditions are necessary to have equation solutions.
The column is initially empty for all protein components, except for the stationary phase concentration of salt that was set to the ionic capacity such that fulfill the electroneutrality condition.
The \textsc{dpfr} and \textsc{cstr} units are initially empty of all components.

\begin{subequations}\label{eq:IC}
    \begin{align}
        & c_0(t=0, z) = 0 && c_i(t=0, z) = 0 \\
        & c_{p,0}(t=0, r) = 0 && c_{p,i}(t=0, r) = 0 \\
        & q_0(t=0, r) = \Lambda && q_i(t=0, r) = 0 \\
        & \tilde{c}_0(t=0, z) = 0 && \tilde{c}_i(t=0, z) = 0 \\
        & \hat{c}_0(t=0, z) = 0 && \hat{c}_i(t=0, z) = 0
    \end{align}
\end{subequations}

\section{Bayesian inference and MCMC sampling}
Consider column models with parameters, $\eta$; the unknown parameters $\theta = (\eta; \tilde{\sigma}^2) \in \mathbb{R}^n$ are estimated from the experimental data, $y \in \mathbb{R}^{d\times m}$, where $n$ is the parameter dimension, $d$ is the number of observation points and $m$ is the number of components.

\subsection{Bayes theorem}
In the Bayesian framework, inference conclusions are made in terms of probabilities, which is used as the fundamental measurement of uncertainties.
Prior distribution, $p(\theta)$, likelihood, $p(y|\theta)$, marginal distribution, $p(y)$, and posterior distribution, $p(\theta|y)$, are related to each other by the Bayes theorem \citep{gelman2014bayesian}:

\begin{equation}
    p(\theta|y)=\frac{p(y|\theta) p(\theta)}{p(y)} 
    \label{eq:Bayes}
\end{equation}
where $p(y) = \int p(y, \theta)\, \dd \theta  = \int p(\theta) p(y|\theta)\, \dd \theta$.
The integral is very hard and computationally expensive to calculate for multi-dimensional distributions.
Therefore, Bayesian inference is often realized by approximating a unnormalized posterior distribution (Eq.~\eqref{eq:likelihood}) via sampling (\ie, \textsc{mcmc}).
Technically, the \textsc{mcmc} sampling allows to calculate the constant value, $p(y)$.
However, this is pointless, as 1) it is sufficient to sample from Eq.~\eqref{eq:likelihood} and 2) integrating $p(y)$ takes the similar computational effort with calculating the sought posterior distribution, $p(\theta|y)$.

\begin{equation}
    p(\theta|y) \propto p(y|\theta) p(\theta) 
    \label{eq:likelihood}
\end{equation}

The probabilistic model used in the present study, which is encoded in the likelihood function, $p(y|\theta)$, assumes that the measurements, $y^T = \{y_{ij}\}$, are corrupted with Gaussian noise:
\begin{equation}
y \sim \mathcal{N}(\tilde{y}(\eta), \Gamma) \label{eq:bayes_probab_model}
\end{equation}

The expected value of this distribution is given by the solution of the mechanistic model, $\tilde{y}^T(\eta) = \{\tilde{y}_{ij}(\eta)\}$, and parameters $\eta$, which can be a subset of the full parameter vector $\theta$.
For simplicity, it is assumed that the noise on each measurement is known, identical and independent, \ie, $\Gamma = \tilde{\sigma}^2 I$ for a given $\tilde{\sigma} > 0$. 
In this case, the likelihood function is given by Eq.~\eqref{eq:bayes_simplified} with $\theta = \eta$.

\begin{equation}\begin{split}
    p(y|\theta) & = \prod_{j=1}^d p(y_j|\theta)  \\
    & \propto \frac{1}{ \pbk{2 \pi \Gamma}^{md/2} } \exp \pbk{ - \frac{1}{2} \sum_{i=1}^m (y_i - \tilde{y}_i(\eta))^T \Gamma^{-1} (y_i - \tilde{y}_i(\eta)) }
    \label{eq:bayes_simplified}
\end{split}\end{equation}

The noise level, $\tilde{\sigma}^2$, can be estimated from the sample variance at a point, $\eta_0$, in the parameter space:

\begin{equation}
    \tilde{\sigma}^2 \approx \tilde{\sigma}_0^2 = \dfrac{\sum_{i=1}^m \pbk{y_i - \tilde{y}_i (\eta_0)}^T \pbk{ y_i - \tilde{y}_i (\eta_0)} }{md-n} 
    \label{eq:sample_variance}
\end{equation}
However, this approach has a critical disadvantage that all subsequent inference steps depend on $\eta_0$ and that it is difficult to find a suitable point $\eta_0$ as a priori.
The problem is circumvented by considering $\tilde{\sigma}^2$ as an additional parameter to be estimated, \ie, $\theta^T = \left( \eta, \tilde{\sigma}^2 \right)$.
A factorization of the prior $p(\theta) = p(\eta) p(\tilde{\sigma}^2)$ yields the posterior density:

\begin{equation}
p(\theta|y) = p(\eta, \tilde{\sigma}^2|y) \propto p(y|\eta, \tilde{\sigma}^2) \, p(\eta) \, p\left(\tilde{\sigma}^2\right) \label{eq:posterior_sigma_mu}
\end{equation}
The inverse Gamma distribution $\text{IG}(\alpha_0, \beta_0)$ with parameters $\alpha_0, \beta_0 > 0$ is taken as prior $p(\tilde{\sigma}^2)$, Eq.~\eqref{eq:sigma_prior}, because then the (conditional) posterior is again an inverse Gamma distribution (conjugate prior, see \cite{gelman2014bayesian}).

\begin{equation}
p\left(\tilde{\sigma}^2\right) \propto \left( \tilde{\sigma}^2 \right)^{-\alpha_0 - 1} \exp\left( -\frac{\beta_0}{\tilde{\sigma}^2} \right) \label{eq:sigma_prior}
\end{equation}
Using small values for $\alpha_0, \beta_0$, \eg, $0.1$, results in a weakly informative prior.
The prior for the column model parameters, $p(\eta)$, will be specified later.

\subsection{Markov Chain Monte Carlo}
Based on the Markov chain theory to generate chains, \textsc{mcmc} is able to sample from complicated distributions.
\textsc{mcmc} is an indirect sampling method; it generates a sequence of random states, each of which depends on the previous one.
The more states that are collected, the more closely the distribution of the samples matches the target distribution.
Various \textsc{mcmc} algorithms have been developed, which mainly differ in computational complexity, robustness, and speed of convergence.

\subsubsection{Metropolis algorithm}\label{sec:MetropolisAlgo}
Metropolis algorithm \citep{metropolis1953equation} is one of the blocking bricks.
It is a random-walk algorithm with Gaussian proposal for sampling the parameters $\theta$.
The algorithm proceeds as follows:

\begin{enumerate}
    \item Initialize a starting point, $\eta^0$, for example, from the prior distribution; construct a covariance matrix, $\Sigma$, for the proposal distribution (Gaussian distribution in the present work).
    \item For $k = 1, 2, \dots$:
        \begin{itemize}
            \item Based on the previous sample, $\eta^k$, a candidate $\tilde{\eta}$ is drawn from the Gaussian proposal distribution, $\mathcal{N}(\eta^k, \Sigma)$.
            \item A ratio $\gamma(\tilde{\eta}, \eta^k)$ of posterior distributions of the candidate, $\tilde{\eta}$, and the previous sample, $\eta^k$, with respect to the desired target distribution is calculated:
\begin{equation}\begin{split}
        \gamma(\tilde{\eta}, \eta^k)  & = \frac{p(\tilde{\eta}|y, \tilde{\sigma}^2)}{p(\eta^k|y, \tilde{\sigma}^2)} = \frac{p(y|\tilde{\eta}, \tilde{\sigma}^2)\, p(\tilde{\eta})}{p(y|\eta^k, \tilde{\sigma}^2)\, p(\eta^k)} \\
        & = \exp \pbk{ - \frac{1}{2} \pbk{\mathcal{S}_{\tilde{\eta}} - \mathcal{S}_{\eta^k}} } \frac{p(\tilde{\eta})}{p(\eta^k)}
    \label{eq:ratio}
\end{split}\end{equation}
The sum of squares in Eq.~\eqref{eq:ratio}, $\sum_{i=1}^m (y_i - \tilde{y}_i(\eta))^T \Gamma^{-1} (y_i - \tilde{y}_i(\eta))$, is denoted as $\mathcal{S}_\eta$.  

            \item The candidate is conditionally accepted with the following probability, where the random number, $\beta$, is drawn from the uniform distribution on the interval $[0,1]$.
        \begin{equation}
            \eta^{k+1} = \left\{  \begin{array}{l@{\quad}l}
            \tilde{\eta} & \beta^k \leqslant \min\pbk{1, \gamma(\tilde{\eta}, \eta^k)} \\
            \eta^k & \text{otherwise} \end{array} \right.
        \end{equation} 

            \item The index $k$ is increased by one and the procedure is repeated until a stopping criterion is satisfied (\eg, a predefined number of samples is reached).
        \end{itemize}
\end{enumerate}

Eventually, a sequence of random samples whose distribution approximates the target density is obtained.
In implementation, a portion of samples (\eg, \SI{25}{\percent}) could be discarded as \textit{burn-in} to diminish the influence of the starting point, $\eta^0$.
Though the Metropolis algorithm has simple and easy to implement features, it has low efficiency.
In sampling, numerous candidates can be rejected, resulting slow convergence to the target distribution.
Further, when one chain is trapped into a local mode, it might never converge to the target density.
Hence, several enhancements have been proposed in the literature.
In this study, an adaptive Metropolis strategy and a delayed rejection \citep{haario2006dram} are applied to alleviate the drawbacks.

\subsubsection{Gibbs sampler}
The Metropolis algorithm described above only samples parameters of the column models $\eta$, but the parameter set also includes $\tilde{\sigma}^2$, which has not been sampled so far.
Including $\tilde{\sigma}^2$ in the Metropolis algorithm, \ie, sampling from the full parameter set $\theta$, is straightforward.
However, the prior $p(\tilde{\sigma}^2)$ has been explicitly chosen such that the conditional posterior $p(\tilde{\sigma}^2 | \eta, y)$, Eq.~\eqref{eq:sigma_conditional_posterior}, is analytically tractable.

\begin{equation}
\tilde{\sigma}^2| y, \eta \sim \text{IG} \left( \alpha_0 + \frac{md}{2}, \beta_0 + \frac{\mathcal{S}_\eta}{2} \right) \label{eq:sigma_conditional_posterior}
\end{equation}

Since it is possible to directly sample from inverse gamma distributions \citep{Marsaglia2000Gamma}, it is more efficient to embed the Metropolis algorithm sampling $\eta$ into a Gibbs sampler that samples $\tilde{\sigma}^2$ directly, than applying the Metropolis algorithm on the full parameter set $\theta$.

A (blocked) Gibbs sampler iteratively sweeps over subsets of parameters and sequentially updates each parameter subset while fixing the remaining ones.
To this end, each subset is drawn from its respective posterior density conditioned on the remaining parameters.
In the case presented above, the Gibbs sampling algorithm starting with iteration counter $k = 0$ and some initial point $\theta^0$ is given as follows.

\begin{enumerate}
\item Sample $\eta^{k+1}$ from $p(\eta | \left(\tilde{\sigma}^2\right)^{k}, y) \propto p(y | \eta, \tilde{\sigma}^2\!=\!\left(\tilde{\sigma}^2\right)^k) \, p(\eta)$ using one iteration of the Metropolis algorithm from Sec.~\ref{sec:MetropolisAlgo}.

\item Sample $\left(\tilde{\sigma}^2\right)^{k+1}$ from $p\left(\tilde{\sigma}^2 | \eta\!=\!\eta^{k+1}, y\right)$ by drawing from the inverse Gamma distribution $\text{IG} \left( \alpha_0 + \frac{md}{2}, \beta_0 + \frac{\mathcal{S}_{\eta^{k+1}}}{2} \right)$. 

\item Increase counter $k$ by one and repeat the procedure from the first step until a stopping criterion is satisfied.
\end{enumerate}

\subsubsection{Adaptive Metropolis and delayed rejection strategies}
The convergence of the Metropolis algorithm can be accelerated by adapting shape of the proposal distribution (\ie, Gaussian distribution in this study), as determined by the covariance matrix, $\Sigma$.
A typical choice for generating the initial variance-covariance matrix, $\Sigma_0$, is Fisher information matrix.
Another practical option is to run a pre-simulation beforehand and then approximate $\Sigma_0$ from the collected samples.
After running sampling at a fixed length, the covariance matrix $\Sigma$ is then adapted in regular intervals, based on the history of the Markov chain.
Further details can be found in \cite{he2019bayesian}.

The sampling robustness and efficiency can be enhanced by delaying the rejection of candidates.
Instead of discarding a proposal, a next stage of the Metropolis algorithm is performed with a shrunken covariance matrix, $a \Sigma$.
Accepting more candidates with cautious moves have benefits of adapting better shape of $\Sigma$ matrix at the very beginning, and evolve the proposal covariance matrix faster to the target density.
Further details can be found in \citep{haario2006dram}.

\section{Clustering-based $\widehat{R}$ diagnostic}

\subsection{$\widehat{R}$ diagnostic}
Assuming that collected samples from multiple chains here render an unimodal target so far.
Consider that the samples collected from $p$ multiple chains are denoted as $\Phi \in \mathbb{R}^{k\times n\times p}$; samples for an estimated parameter $\theta_\ell\, (\ell \in \{1,\dots,n\})$ are labelled as $\Phi^\ell_{\varkappa, \varrho}, (\varkappa\in\{1,\dots,k\}, \varrho \in\{1, \dots, p\}$).
As introduced, the potential scale reduction factor \citep{gelman2014bayesian} is applied in this study to assess convergence conditions.
It is a square root of the ratio of sample variances, 
which is based on between- and within-chain variances:
\begin{subequations}\begin{align}
    \widehat{R}_\ell & = \sqrt{\frac{\widehat{\text{var}}^{+}(\mathcal{W}, \mathcal{B})}{\mathcal{W}}} \label{eq:gelman} \\
    \mathcal{B} & = \frac{k}{p-1} \sum_{\varrho=1}^p \left( \bar{\Phi}_{.\varrho}^\ell - \bar{\Phi}_{..}^\ell \right)^2 \label{eq:B} \\
    \mathcal{W} & = \frac{1}{p} \sum_{\varrho=1}^p s_\varrho^2 \label{eq:W}
\end{align}\end{subequations}
Note that the \emph{burn-in} part has been discarded in length of the chains $k$ here.
$\bar{\Phi}_{.\varrho}^\ell = \frac{1}{k} \sum_{\varkappa=1}^k \Phi_{\varkappa, \varrho}^\ell$ is the mean value of the chain $\varrho$.
$\bar{\Phi}_{..}^\ell$ is the mean value of the mean vector of $p$ chains, $\frac{1}{p} \sum_{\varrho=1}^p \bar{\Phi}_{.\varrho}^\ell$.
Thus, $\mathcal{B}$ defines the between-chain variance; $s_\varrho^2$ denotes the within-chain variance of the chain $\varrho$.
The sample variance, $\widehat{\text{var}}^{+}(\mathcal{W}, \mathcal{B})$, is estimated by a weighted average of $\mathcal{W}$ and $\mathcal{B}$, namely
\begin{equation}
    \widehat{\text{var}}^{+}(\mathcal{W}, \mathcal{B}) = \frac{k-1}{k} \mathcal{W} + \frac{1}{k} \mathcal{B}
\end{equation}

Upon convergence of the \textsc{mcmc} algorithm diagnosed by $\widehat{R}$ diagnostic, the samples collected, since then, from the multiple chains can be mixed up to approximate the target distribution.
The effective number of independent draws (\ie, \emph{effective sample size}) for any optimized parameter $\theta_\ell$ can be estimated from Eq.~\eqref{eq:eff}.
\begin{equation}
    n_{\text{eff}} = \frac{p\, k}{1 + 2 \displaystyle\sum_{t=1}^{\infty} \rho_t}
    \label{eq:eff}
\end{equation}
where $\rho_t$ is the autocorrelation of the mixed-up chain of the parameter $\theta_\ell$ at lag $t$; 
In practice, however, we barely have a finite simulation length, so the calculation has to be approximated.

\subsection{$K$-means clustering}
Clustering is a powerful tool for analyzing structures of datasets.
Various clustering algorithms are available, all of which take slightly different approaches and might produce slightly different results.
Centroid-based clustering, $K$-means, and its variants are used in this study.
We firstly introduce the idea with the demonstration of $K$-means.
The $K$-means algorithm here can be replaced with its variants and even hierarchical methods without losing generality.
$K$-means aims at partitioning observations into $K$ clusters, in which each observation belongs to the cluster with the nearest mean.
In this study, observations are $p$ multiple chains from the \textsc{mcmc} sampling.

The basic idea is that when the posterior distribution (cf.~Eq.~\eqref{eq:likelihood}) is rough and multi-modal, we could partition it into $K$ sub-posterior distributions:

\begin{equation}
    p(\theta|y) = \prod_{j=1}^K \omega_j\, p_j(\theta|y)
    \label{eq:subposterior}
\end{equation}
where $\omega_j$ are the weights.
For each sub-posterior distribution, at least one of $p$ chains are triggered to sample its target density.
Therefore, the $K$ clusters of chains definitely have different convergence rate; different cluster of chains might be transitional between each other with rather low probability.
However, the \textsc{mcmc} sampling can only be stopped (\eg, $\widehat{R}$ diagnostic) when each chain traverse all the target densities of the posterior distribution.

The key idea is that before applying $\widehat{R}$ diagnostic to assess convergence of all the chains, $K$-means clustering \citep{lloyd1982least} is firstly applied to partition the $p$ chains into $K$ groups. 
For each cluster of $K$ (which is assumed to be unimodal), the collected samples from \textsc{mcmc} are assessed with $\widehat{R}$ diagnostic.
It is asserted to be converged when $\widehat{R}$ is below a threshold (\eg, 1.10) for each of $K$ cluster, instead waiting extreme long time to let each chain traverse all the modes.
An extension of the idea to $K$-means variants and even hierarchical methods is straightforward.

Assume the dataset requires classification is denoted as $\mathcal{X} \in \mathbb{R}^{p\times q}$, that is consisted of $p$ chains, each chain $x_\varrho$ has dimension of $q$, where $\varrho \in \{1, \dots, p\}$:
\begin{equation}
    \mathcal{X} = \sbk{x_1, \dots, x_\varrho, \dots, x_p}^T
    \label{eq:kmeansX}
\end{equation}
The centroid of each cluster is denoted as $\mu_j, j \in \{1, \dots, K\}$, such that $\mathcal{U} \in \mathbb{R}^{K\times q}$:
\begin{equation}
    \mathcal{U} = \sbk{\mu_1, \dots, \mu_j, \dots, \mu_K}^T
\end{equation}
Binary indicators $z_\varrho^j$ associated to $x_\varrho$ are introduced such that $z_\varrho^j = 1$ if $x_\varrho$ belongs to the cluster $j$, $z_\varrho^j = 0$ if not.
The Euclidean distance $\mathcal{L}$ between each chain $x_\varrho$ and its belonging centroid $\mu_j$ (a.~k.~a.~distortion) is calculated as:
\begin{equation}
    \mathcal{L}(\mu, z) = \sum_{\varrho=1}^p \sum_{j=1}^K z_\varrho^j\, \norm{x_\varrho - \mu_j}^2
    \label{eq:distortion}
\end{equation}

The $K$-means algorithm is proceeded as follows:
\begin{enumerate}
    \item apply traditional $\widehat{R}_\ell$ to each model parameter $\theta_\ell, \forall \ell \in \{1,\dots,n\}$. If $\widehat{R}_\ell < 1.10\quad \forall \ell$, the \textsc{mcmc} sampling is stalled.
    \item otherwise
        \begin{enumerate}
        \item specify a value of $K$. Choosing a proper $K$ depends on datasets and is primarily subjective. Finding the optimal $K$ value shall be introduced in Section \ref{sec:optimalK}.
        \item randomly generate $K$ centroids in $q$ dimension, $\mathcal{U}$. The randomly generated centroids can be arbitrarily bad. $K$-means++ algorithm is used to address this problem. A procedure is implemented to initialize the cluster centers before proceeding with the $K$-means iterations. With the $K$-means++ initialization, it is guaranteed to find a solution that is $\mathcal{O}(\log K)$ competitive to the optima \citep{arthur2006k}.
        \item assign chains $\mathcal{X}$ to the generated centroids $\mathcal{U}$, by calculating Euclidean distance, Eq.~\eqref{eq:distortion}, from every chain to every centroid. We then take the closet centroid for each chain and assign it to that cluster.
        \item re-calculate the centroid $\mathcal{U}$. Thereafter, we compute the mean for each assigned cluster and place the centroid on the mean value.
        \item repeat step (c) and (d) until convergence.
    \end{enumerate}
\end{enumerate}

One thing deferred here is what actually $\mathcal{X} \in \mathbb{R}^{p\times q}$ is?
As the chain matrix collected from sampling, $\Phi^{k\times n\times p}$, is multi-dimensional, which can not be transferred to use $K$-means algorithms directly.
A mapping function is applied to reduce dimension such that to have $\mathcal{X}$ in Eq.~\eqref{eq:kmeansX}:
\begin{equation}
    f: \pbk{ \Phi \in \mathbb{R}^{k\times n\times p} } \mapsto \pbk{ \mathcal{X} \in \mathbb{R}^{p\times q} }
    \label{eq:mapping}
\end{equation}
The samples for each model parameter, $\theta_\ell^\varrho, \ell \in \{1, \dots, n\}$, of chain $\varrho$ with length of $k$ are represented with probability density, $h_\ell^\varrho$.
Kernel smoothing function estimate can be used to generate discrete points of $h_\ell^\varrho \in \mathbb{R}^q$, where $q = n\, s$ and $s$ is an user-defined parameter:
\begin{equation}
    h_\ell^\varrho = f \pbk{\Phi_{. \varrho}^\ell, s}
    \label{eq:ksd}
\end{equation}
Then, the kernel estimates of $n$ parameters are concatenated into a row vector, 
\begin{equation}
    x_\varrho = [h^\varrho_1, \dots, h^\varrho_\ell, \dots, h^\varrho_n]
\end{equation}
Lastly, the dataset for classification, $\mathcal{X} \in \mathbb{R}^{p\times q}$, in constructed from 
\begin{equation}
    \mathcal{X} = \sbk{x_1, \dots, x_\varrho, \dots, x_p}^T
\end{equation}
In the framework of Bayesian inference, kernel density might be a good choice with high sensitivity for constructing dataset for mode clustering.

\subsection{Optimal $K$}\label{sec:optimalK}
The choice of $K$ is not universal.
As observed, when $K$ increases the distortion $\mathcal{L}$ decreases, until it reaches $0$ when $K = p$ (each data is the center of its own).
Therefore, a penalty item over $K$ can be added to the original distortion, $\mathcal{L}$, resulting in a new minimization problem:
\begin{equation}
    \mathcal{L}(\mu, z, K) = \sum_{\varrho=1}^p \sum_{j=1}^K \norm{x_\varrho - \mu_j}^2 + \lambda K
\end{equation}
However, the choice of $\lambda$ is also arbitrary or heuristic again.
Another option is that running the $K$-means algorithm over a linear range of $K$ values, then an elbow method is used to determine the optimal $K$ value from the $K \text{vs.}~\mathcal{L}$ figure.
The latter scheme is adopted in our study.

\section{Case study}
Fig.~\ref{fig:model_scheme} depicts a general picture of column models; models can be omitted to be suitable for different case studies in the chromatographic field. For instance, the \textsc{grm} can be omitted to account for column bypass experiments to identify hold-up volumes. 
Description of guard columns before/after the main column are adequate with \textsc{dpfr} in some scenarios.

The multi-modal identification using Bayesian inference and \textsc{mcmc} sampling and clustering-based $\widehat{R}$ diagnostic is demonstrated in ion-exchange chromatography with lysozyme as target protein.
The \textsc{sma} model is required for describing gradient elution experiments and, consequently, the \textsc{grm} also includes a salt component. 
Surface diffusion in Eq.~\eqref{eq:GRM} is not considered in the porous volume.

The unknown model parameters are estimated using Bayesian inference.
Parameters of the column models $L^\textsc{dpfr}$, $D^\textsc{dpfr}_\text{ax}$, $V^\textsc{cstr}$, $\varepsilon_c$, $D_\text{ax}$, $\varepsilon_p$, $D_p$, $k_f$, $k_\text{eq}$, $\nu$, $\sigma$ are estimated.
The remaining parameters, namely column length, $L$, and particle diameter, $r_p$, are known.
The linear flow rates in the \textsc{dpfr}, $u^\textsc{dpfr}$, and column, $u_\text{int}$, are calculated from the volumetric flow rate and the respective cross-sectional areas.
The time-dependent inlet concentrations of salt and protein, $c^{\text{in}}_i$, are used as boundary conditions.
The outlet of previous unit is used as the inlet of the downstreaming unit.
Uniform prior distributions between upper and lower boundaries are defined for each parameter.
The parameters of the inverse Gamma distribution in Eq.~\eqref{eq:sigma_prior} are chosen as $\alpha_0 = 0.5$ and $\beta_0 =  0.5\, \tilde{\sigma}_0^2$.
Posterior distributions are visualized by applying a kernel density estimator, \textsc{matlab} command $\mathtt{ksdensity}$, as well as the mapping of Eq.~\eqref{eq:ksd}.

The orders of magnitude of the model parameters are harmonized by applying a logarithmic transformation previous to parameter estimation.
Consequently, the estimation procedure operates on the exponents, $\rho = \ln(\eta)$, instead of the original parameter values, $\eta = \exp(\rho)$.
This changes the posterior density due to the change of variables theorem.

\begin{equation}
p( \rho, \tilde{\sigma} |y) \propto p(y|\exp(\rho), \tilde{\sigma}) \, p(\exp(\rho)) \, p\left(\tilde{\sigma}^2\right) \, \prod_{i=1}^n \exp(\rho_i)  \label{eq:posterior_transformed}
\end{equation}
Note that $\exp(\rho)$ is understood component-wise in Eq.~\eqref{eq:posterior_transformed}.
The inverse transform, $\eta = \exp(\rho)$, is applied after the estimation procedure.

The \textsc{grm} is solved using the chromatography analysis and design toolkit (\textsc{cadet}) which is freely available as open-source software (\url{https://github.com/modsim/CADET}).
The column is discretized using \num{100} axial elements, and the particles are discretized using \num{10} radial elements.
Parameter estimation using Bayesian inference and \textsc{mcmc} sampling is implemented in a separate package that is available under the same conditions (\url{https://github.com/modsim/CADET-MCMC}).

\section{Experiment}
Prepacked \SI{1}{\milli\liter} \textsc{sp} Sepharose \textsc{ff} HiTrap columns from \textsc{ge} Healthcare were used.
They have a length of $L = \SI{2.5e-2}{\meter}$ and an inner diameter of $d_c = \SI{7e-3}{\meter}$.
The particle radius is $r_p = \SI{4.5e-5}{\meter}$.
Experiments were performed using an \textsc{\"akta} pure \SI{25}{\liter} system (\textsc{ge} Healthcare, Little Chalfont, \textsc{uk}) controlled by the Unicorn 6.4 software (\textsc{ge} Healthcare).
The flow rate of $Q = \SI{8.33e-9}{\cubic\metre\per\second}$ was kept constant in all experiments.
Ultrapure water was prepared with an arium pro \textsc{vf} system (Sartorius, G\"ottingen, Germany) and used in all experiments.
Monosodium phosphate dihydrate, and disodium phosphate dihydrate (buffer) were purchased from Carl Roth (Karlsruhe, Germany).
Lysozyme from chicken egg white (L4919) was purchased from Sigma Aldrich (St.~Louis, \textsc{usa}).

In the elution experiments, \SI{0.1}{\milli\liter} of \SI{0.2}{\milli\Molar} lysozyme solution, prepared by dissolving \SI{0.2}{\milli\mole} lysozyme in \SI{20}{\milli\Molar} phosphate buffer with a pH of 7.0, was loaded to the column.
The column was then equilibrated with \SI{20}{\milli\Molar} phosphate buffer at pH 7.0 for 10 \textsc{cv} (\SI{1155}{\second}).
The elution buffer contained \SI{500}{\milli\Molar} sodium chloride dissolved in \SI{20}{\milli\Molar} phosphate buffer with a final pH of 7.0.

Experiments with lysozyme were performed in a step elution mode. 
First, $c^\text{in}_\text{lyz} = \SI{0.2}{\mole\per\cubic\metre}$ of protein was loaded to the column for $t_\text{load} = \SI{12.6}{\second}$ at a salt concentration of $c^\text{in}_0 = \SI{20}{\mole\per\cubic\metre}$.
The column was then washed for $t_\text{wash} = \SI{578.4}{\second}$ at the same salt concentration. 
The bound proteins were then eluted by increasing the salt concentration from $c^\text{in}_0 = \SI{20}{\mole\per\cubic\metre}$ to $\SI{526}{\mole\per\cubic\metre}$ over time of 10 \textsc{cv} (\SI{1155}{\second}).

\section{Results and discussion}
The lysozyme step elution data are used for estimating the model parameters.
Uniform distributions with upper and lower limits shown in Tab.~\ref{tab:prior_bounds} were used for these parameters.
The ionic capacity of the \textsc{sma}, $\Lambda$, has a distribution as it depends on the total porosity, $\varepsilon_t = \varepsilon_c + (1 - \varepsilon_c) \varepsilon_p$. 

\begin{equation}
\Lambda = \frac{c_{\text{in},s} V^\textsc{tit}}{(1-\varepsilon_t) V^\textsc{col}}\label{eq:lambda}
\end{equation}
In Eq.~\eqref{eq:lambda}, $c_{\text{in},s} = \SI{10}{\mole\per\cubic\metre}$ and $V^\textsc{tit} = \SI{19.25e-6}{\cubic\metre}$ are the titration concentration and volume, and $V^\textsc{col} = \SI{9.62e-7}{\cubic\metre}$ is the volume of the empty column without particles \citep{osberghaus2012determination}.
The ionic capacity was determined using the method described by \cite{huuk2014model}.
\begin{table}[!ht]
    \footnotesize
    \centering
    \caption{Upper and lower bounds of uniform prior distributions.}
    \label{tab:prior_bounds}
    \begin{tabular}{lccccccccccc}
        \toprule
            & $L^\textsc{pfr}$ & $D_\text{ax}^\textsc{pfr}$ & $V^\textsc{str}$ & $D_\text{ax}$ & $\varepsilon_c$ & $k_f$ & $D_p$ & $\varepsilon_p$ & $k_\text{eq}$ & $\nu$ & $\sigma$	\\
        \midrule
        $\min$  & 0.01  & \num{1e-12}   & \num{1e-11}   & \num{1e-12}   & 0.01  & \num{1e-9}    & \num{1e-12}   & 0.40  & \num{1e-3}    & 1.0   & 4.0   \\
        $\max$  & 1.00  & \num{1e-8}    & \num{1e-7}    & \num{1e-7}    & 0.90  & \num{1e-4}    & \num{1e-7}    & 0.99  & 1.0           & 30    & 80	\\
        \bottomrule
    \end{tabular}
\end{table}

\subsection{Clustering and diagnostic of MCMC chains}

\begin{figure}[ht]
    \centering
    \includegraphics[width=0.55\textwidth]{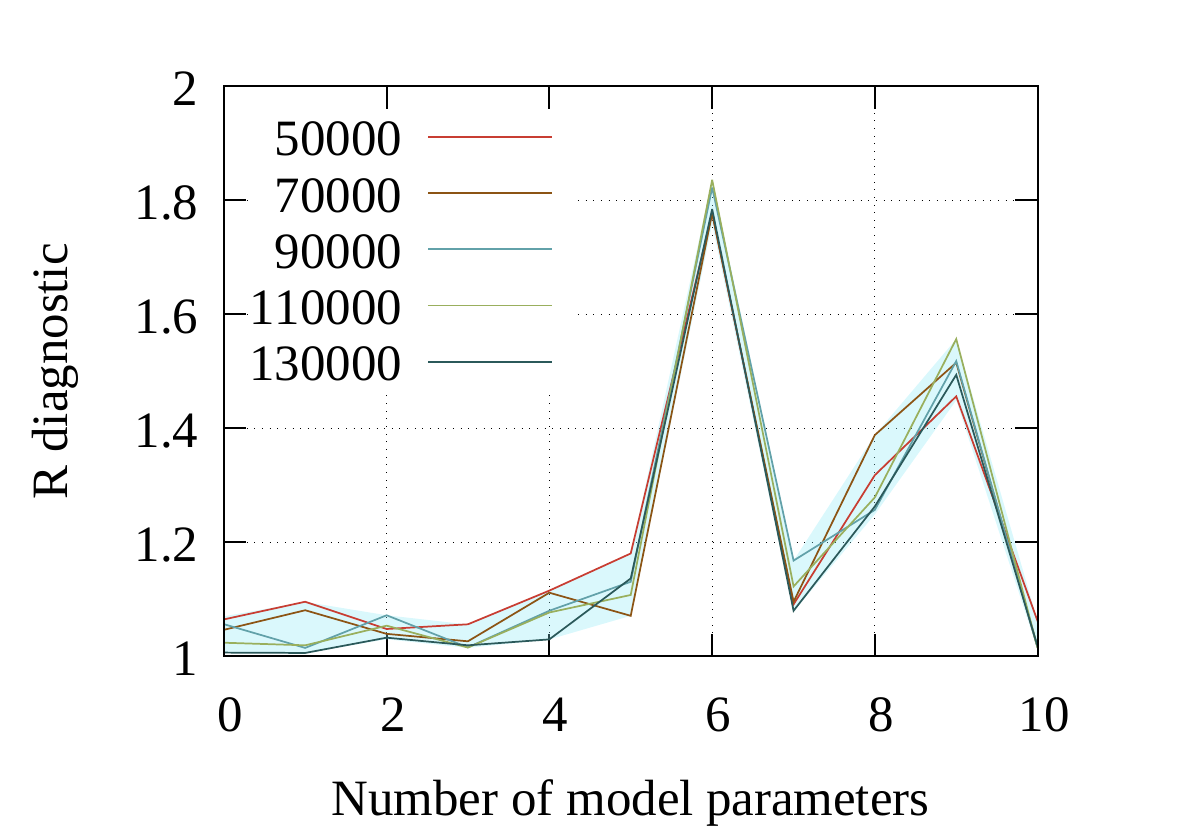}
    \caption{$\widehat{R}$ evolution from sampling length of 50000 to 70000. $0:10$ corresponds to the index of model parameters $L^\textsc{pfr}$, $D_\text{ax}^\textsc{pfr}$, $V^\textsc{str}$, $\varepsilon_c$, $D_\text{ax}$, $\varepsilon_p$, $D_p$, $k_f$, $k_\text{eq}$, $\nu$, and $\sigma$. $\widehat{R}$ values of $D_p$, $k_\text{eq}$ and $\nu$ parameters are higher than \num{1.20}. }
    \label{fig:Revo}
\end{figure}

Six \textsc{mcmc} sampling chains were triggered in this study as we have to consider the capacity of our computing node.
It could apply to more parallel chains without losing generality; it is even more powerful to implement the proposed method on \textsc{gpu}.
Hereafter, the six chains are denoted with randomly generated alphabets for convenience, \ie, $\texttt{e}, \texttt{k}, \texttt{r}, \texttt{u}, \texttt{w}, \texttt{x}$. 
The sampling length of $\texttt{e}, \texttt{k}, \texttt{r}, \texttt{u}, \texttt{w}, \texttt{x}$ chains is $[\num{148582}, \num{143388}, \num{148186}, \num{147701}, \num{147594}, \num{144960}]$, respectively; the total mixing length is $\num{880411}$.
If the $\widehat{R}$ diagnostic is applied to the six parallel sampling chains, it returns a vector with respect to the model parameters:
\begin{equation}
    \widehat{R} = \sbk{1.01, 1.01, 1.03, 1.02, 1.02, 1.13, 1.76, 1.08, 1.28, 1.49, 1.02}
\end{equation}
which definitely turns out to be non-converged, when the threshold is 1.20. 
Interestingly, as observed in Fig.~\ref{fig:Revo}, the $\widehat{R}$ vector is almost unchanged since the sampling length, $\num{50000}$, instead it is even getting worse when the sampling length is longer.
We could bravely speculate that the $\widehat{R}$ vector would not be improved so much even the sampling length is doubled (then, the mixing length is around \num{2e6}).
This is an evidence that the traditional $\widehat{R}$ diagnostic do not behave well when the target distribution is rough and multi-modal.
It might work when the simulation is infinitely long; but, the computation cost is unbearable.
Further, $\widehat{R}$ values of $D_p$, $k_\text{eq}$ and $\nu$ parameters are higher than \num{1.20}; more attention should be paid to these parameters.

After applying the presented $K$-means based $\widehat{R}$ diagnostic to the respectively collected samples of the $\texttt{e}, \texttt{k}, \texttt{r}, \texttt{u}, \texttt{w}, \texttt{x}$ chains, it indicates that there are $K=3$ clusters of the six chains (cf.~\ref{sec:appd} for finding the optimal $K$ value).
To be specific, $\cbk{\texttt{w}}$ is a group that differs with all the remaining chains; 
$\cbk{\texttt{e}, \texttt{r}, \texttt{x}}$ is a group that renders a mode, while $\cbk{\texttt{k}, \texttt{u}}$ renders the another one.
The $\widehat{R}$ diagnostic vectors for $\cbk{\texttt{e}, \texttt{r}, \texttt{x}}$ and $\cbk{\texttt{k}, \texttt{u}}$ groups are respectively, 
\begin{equation}
    \begin{bmatrix}
        1.01, 1.01, 1.01, 1.00, 1.02, 1.03, 1.01, 1.01, 1.01, 1.03, 1.01 \\
        1.02, 1.00, 1.00, 1.00, 1.05, 1.00, 1.18, 1.06, 1.23, 1.19, 1.01
    \end{bmatrix}
\end{equation}
Note that the $\widehat{R}$ value of $\cbk{\texttt{w}}$ group can not be calculated as there is only one chain.
The $\cbk{\texttt{e}, \texttt{r}, \texttt{x}}$ turns out to be converged as the highest value is below \num{1.10}.
Assuming the clustering of multiple chains is a \emph{priori} and their convergence rates are monitored separately at the very beginning.
When the sampling length is $k = \num{18491}$, $\widehat{R}^{\cbk{\texttt{e}, \texttt{r}, \texttt{x}}} < 1.20$; when $k = \num{39068}$, $\widehat{R}^{\cbk{\texttt{e}, \texttt{r}, \texttt{x}}} < 1.10$. 
Both of them are way shorter than the current sampling length ca.~\num{1.45e5}.

Other variants of the $K$-means method have been cross-validated to see if other $K$ values would be reported.
$K$-medoids is a clustering algorithm which is more robust and resilient to outliers compared to the $K$-means, as it chooses the actual data points as the prototypes.
A subroutine, \texttt{kmedoids}, from the \textsc{matlab} has been used, where the partitioning around medoids (\textsc{pam}) algorithm and the squared Euclidean distance were used.
The $K$ vs.~$\mathcal{L}$ plot was shown in \ref{sec:appd}.
All the results from the $K$-means, $K$-medoids, $K$-medians and dendrogram agree with each other.

\begin{figure}[h]
    \centering
    \includegraphics[width=0.45\textwidth]{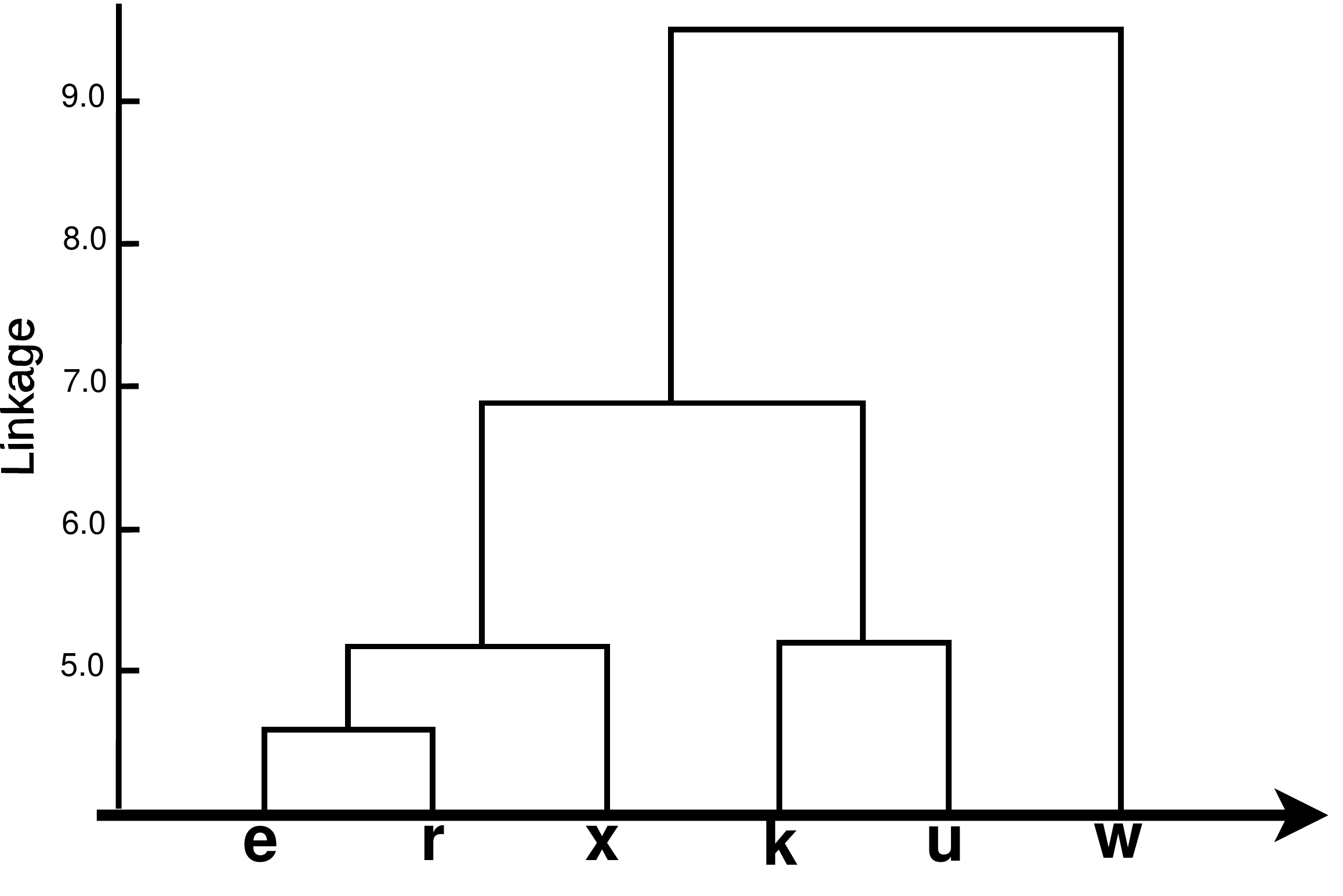}
    \caption{Dendrogram of the sampling chains}
    \label{fig:dendrogram}
\end{figure}
Finding the optimal $K$ value would be a critical point.
The elbow finding algorithm used in this study is not general enough.
However, it does work in our case study.
Optimal $K$ searching strategy would be improved later in our work.
\emph{Dendrogram}, which depicts the hierarchical relationship among objects (\eg, chains) in tree branches, can also be applied here, see Fig.~\ref{fig:dendrogram}.
In this representation, we do need to worry about the determination of $K$ value.
It further elucidates that in the $\cbk{\texttt{e}, \texttt{r}, \texttt{x}}$ group, $\texttt{e}$ and $\texttt{r}$ chains are in closer relation.
The application of the proposed convergence diagnostic can be applied either online or offline.

\subsection{Multi-modal posterior distributions}
The $K=3$ modes of posterior distributions from the $\cbk{\texttt{e}, \texttt{r}, \texttt{x}}$, $\cbk{\texttt{k}, \texttt{u}}$ and $\texttt{w}$ groups are overlapped for comparison here.
For bounded and tailing target distributions, \textsc{mcmc} sampling takes extremely long time to converge to satisfy the $\widehat{R} \approx 1$ criterion.
Thus, for target distributions that are far from normal in high dimensional space, we have generally been satisfied with setting 1.1 or even 1.2 as the threshold.
Mixing length of \num{83167} has been collected from the $\cbk{\texttt{e}, \texttt{r}, \texttt{x}}$ cluster to generate the posterior distributions, while length of \num{39506} collected from the $\cbk{\texttt{k}, \texttt{u}}$ cluster.
Regarding the $\cbk{\texttt{w}}$ cluster where the $\widehat{R}$ diagnostic can not be applied, the chain with \emph{burned-in} length of \num{50000} was used.
The model parameters are catalogued into \emph{system}, \emph{column and particle} and \emph{binding} properties.

\subsubsection{System properties}
\begin{figure}[ht]
    \centering
    \begin{subfigure}{0.33\textwidth}{\includegraphics[width=\textwidth]{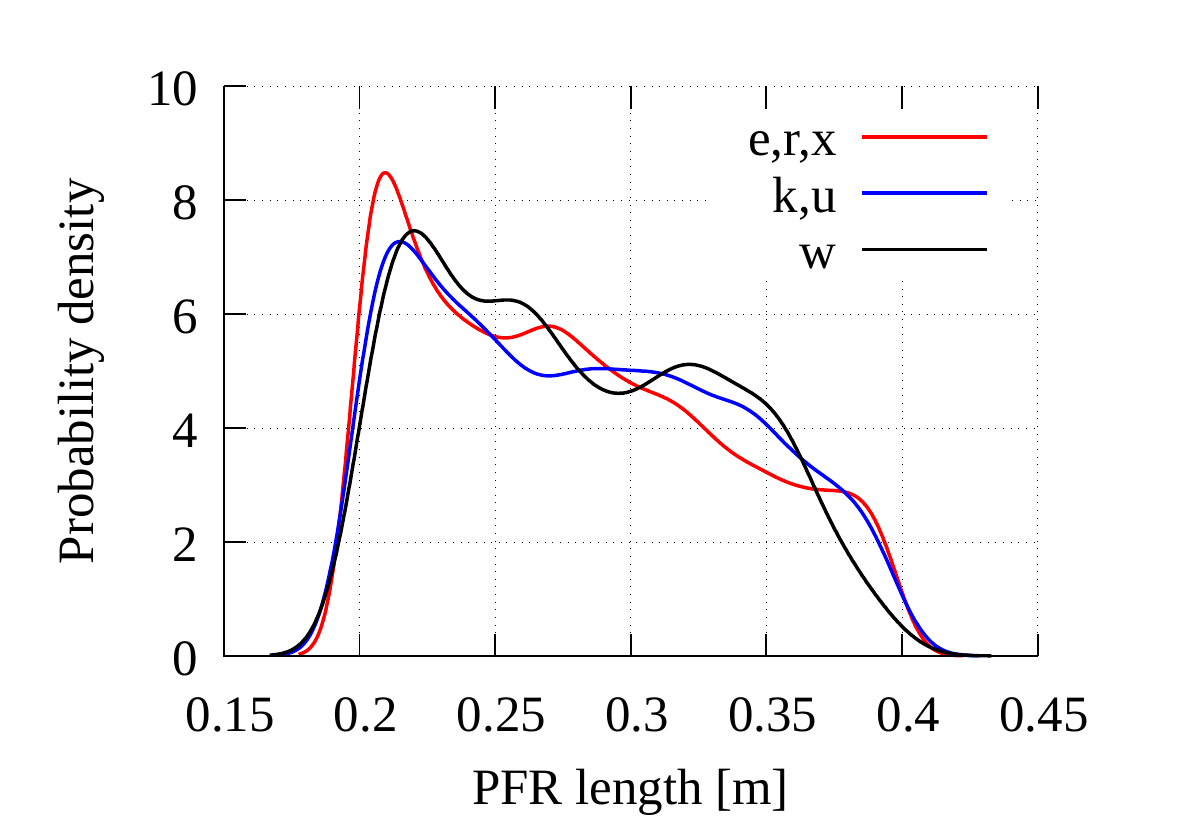}} \caption{\textsc{pfr} length}	\label{fig:postdist_sys_a} \end{subfigure}
    \begin{subfigure}{0.33\textwidth}{\includegraphics[width=\textwidth]{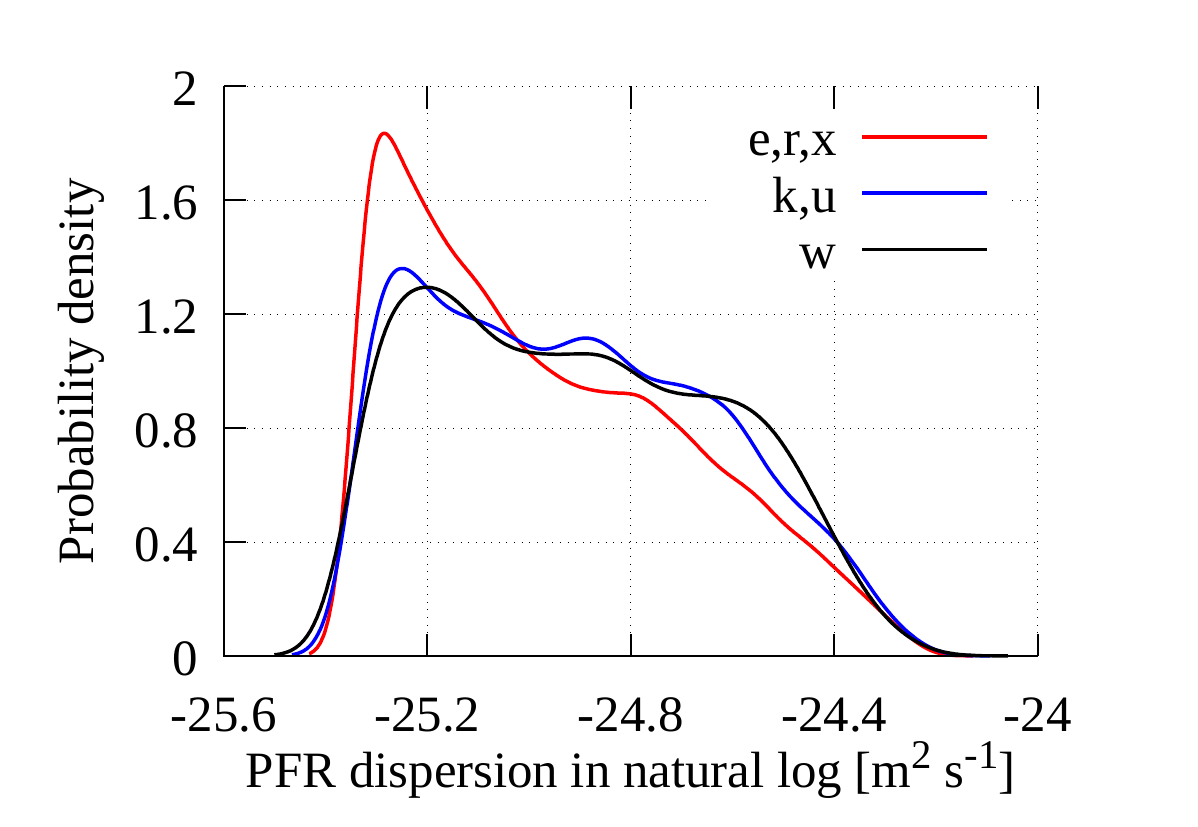}} \caption{\textsc{pfr} dispersion (log)} \label{fig:postdist_sys_b} \end{subfigure}
    \begin{subfigure}{0.33\textwidth}{\includegraphics[width=\textwidth]{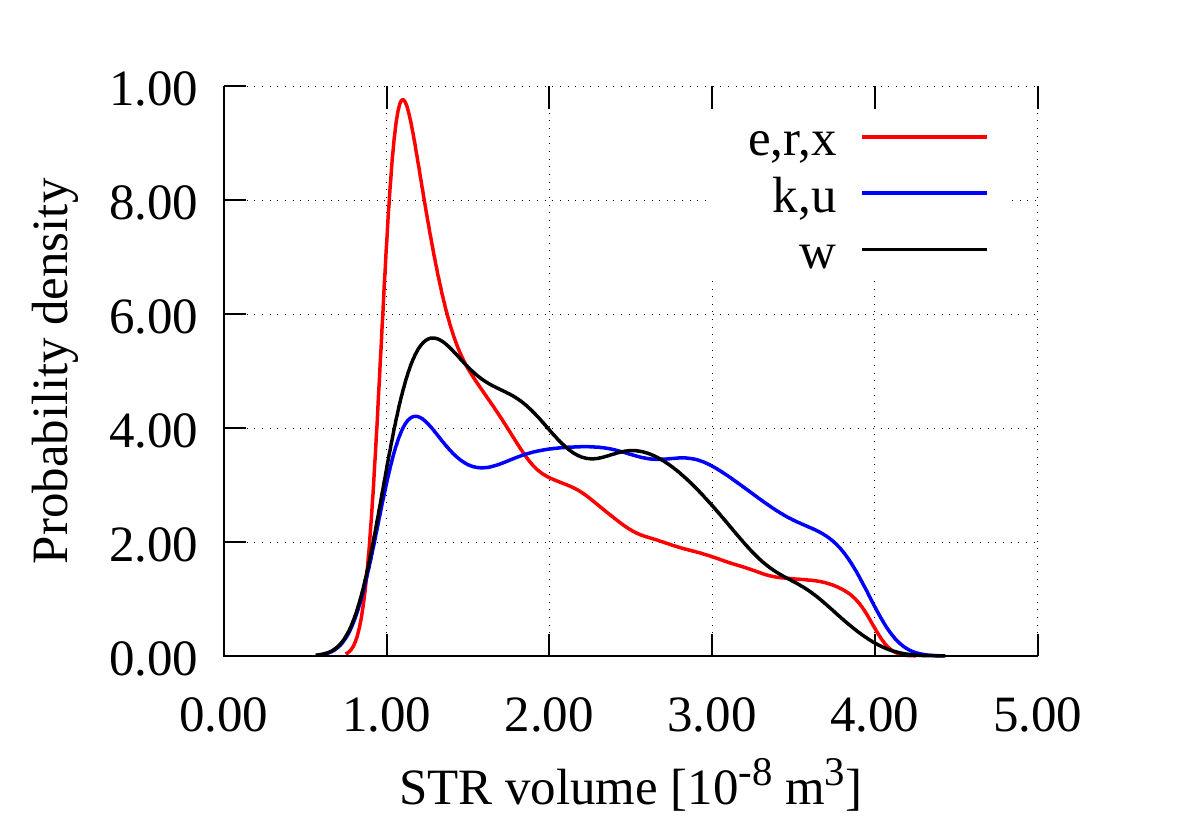}} \caption{\textsc{str} volume} \label{fig:postdist_sys_c} \end{subfigure} 
    \begin{subfigure}{0.33\textwidth}{\includegraphics[width=\textwidth]{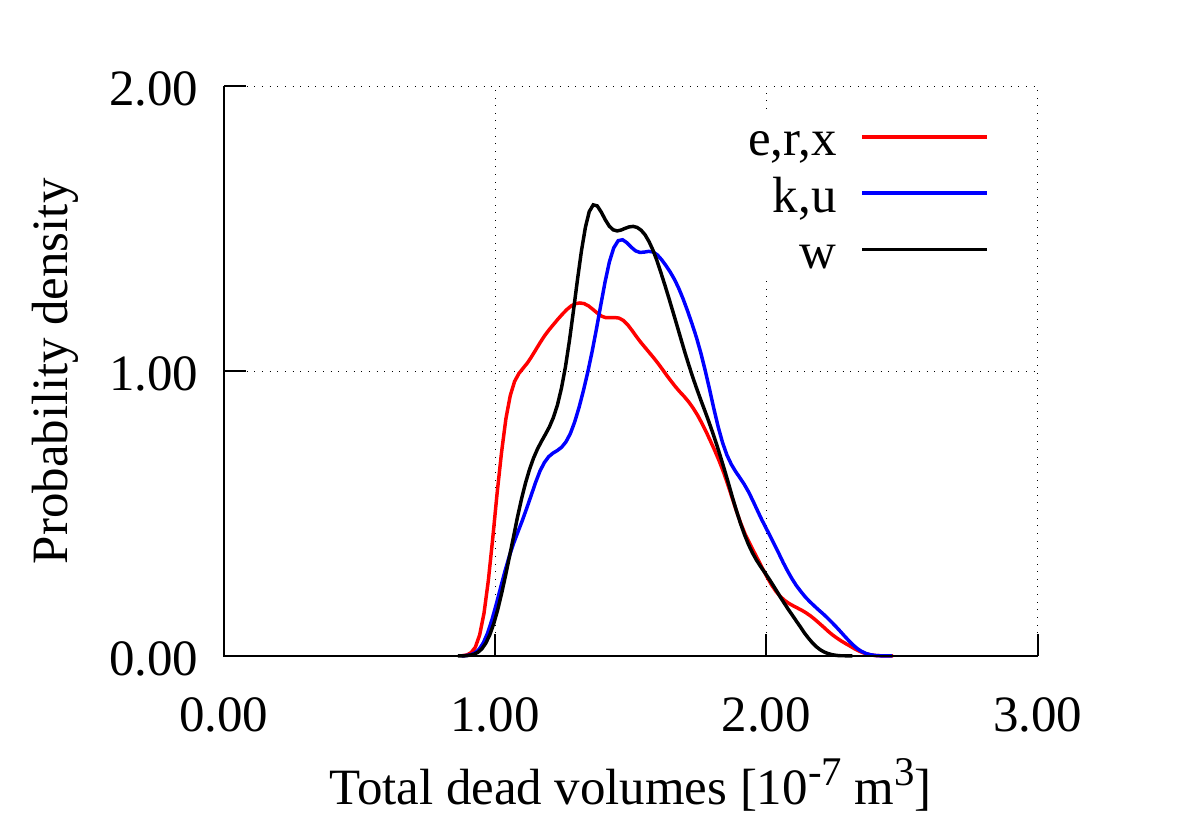}} \caption{Total hold-up volumes} \label{fig:postdist_sys_d} \end{subfigure}
    \caption{} 
    \label{fig:postdist_sys}
\end{figure}

The effects of system volumes within frits, tubing, pumps, valves and detector is characterized by the \textsc{dpfr} and \textsc{cstr} models, which can also be experimentally described with detached column and acetone and dextran as tracers.
However, due to high molecular weight of the dextran, unexpected diffusion behaviour of dextran molecules might be observed.
\cite{karlsson2004aspects} and \cite{persson2004calibration} have reported that the dextran can partially penetrate the particle pores.
Two identical \textsc{dpfr} and \textsc{cstr} units in a mirrored way (cf.~Fig.~\ref{fig:model_scheme}) are implemented in this case study.
For the given flow rate, $Q = \SI{8.33e-9}{\cubic\metre\per\second}$, and tubing radius, $r_t = \SI{2.5e-4}{\meter}$, these two models are uniquely determined by three parameters, namely \textsc{dpfr} length, $L^\textsc{dpfr}$, axial dispersion, $D_\text{ax}^\textsc{dpfr}$, and \textsc{cstr} volume, $V^\textsc{cstr}$.

For the system properties, distributions of $K=3$ modes are consistent with each other, see Fig.~\ref{fig:postdist_sys}.
Therefore, system properties do not render multi-modality and can be well-identified.
The distributions of system properties are slightly rough in the shape of curves, see Fig.~\ref{fig:postdist_sys}.
Longer sampling length could polish the roughness of the distributions (\eg, Fig.~\ref{fig:postdist_sys_a}); but, the roughness does not significantly disturb the inference here.
Dispersion of the system volumes is located around \SI{1e-11}{\square\metre\per\second}, which turns out to be smaller than that of the interstitial and porous volumes.
The systematic hold-up volumes are calculated as $V^\text{tot} = \pbk{A L^\textsc{dpfr} + V^\textsc{cstr} }\times 2$, and the distribution is shown in Fig.~\ref{fig:postdist_sys_d}.
As seen in Fig.~\ref{fig:postdist_sys_d}, the total hold-up volumes is between \SI{1e-7}{\cubic\metre} and \SI{2e-7}{\cubic\metre}.

\subsubsection{Column and particle properties}
\begin{figure}[h]
    \centering
    \begin{subfigure}{0.33\textwidth}{\includegraphics[width=\textwidth]{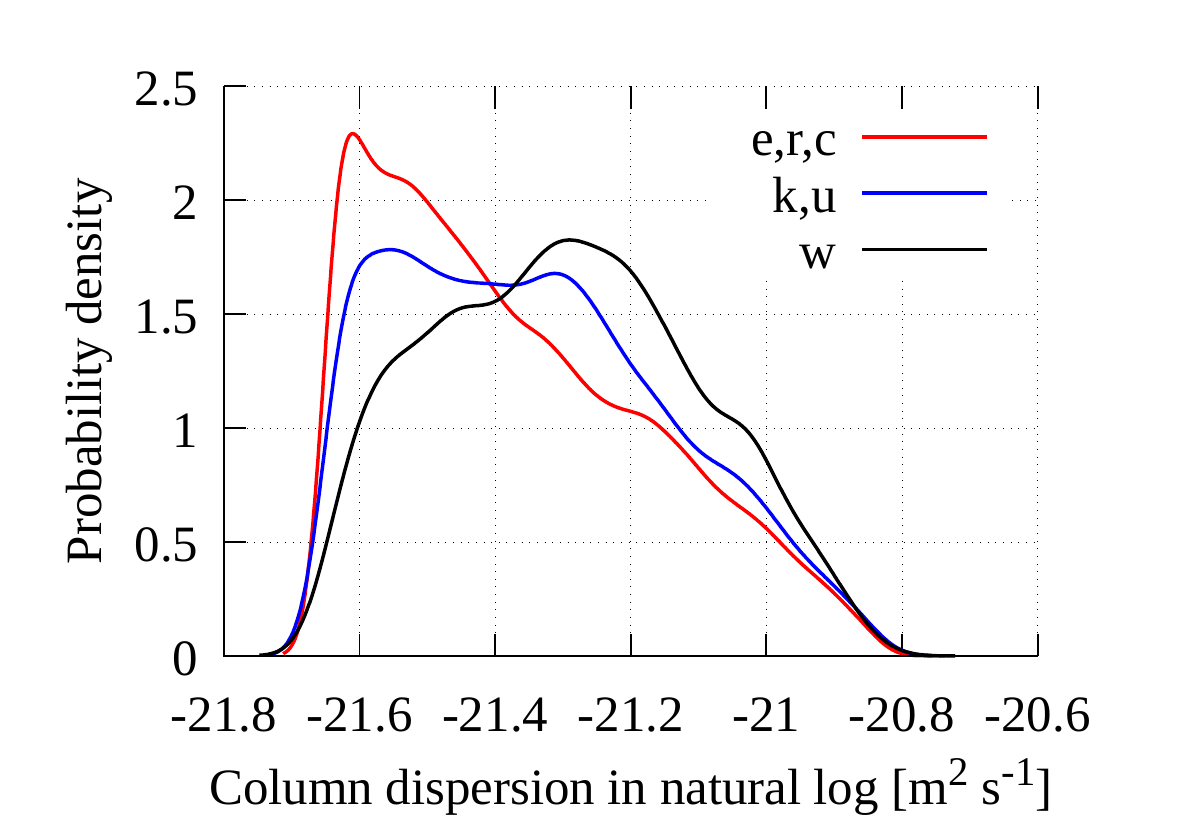}} \caption{Column dispersion (log)} \label{fig:postdist_col_a} \end{subfigure}
    \begin{subfigure}{0.33\textwidth}{\includegraphics[width=\textwidth]{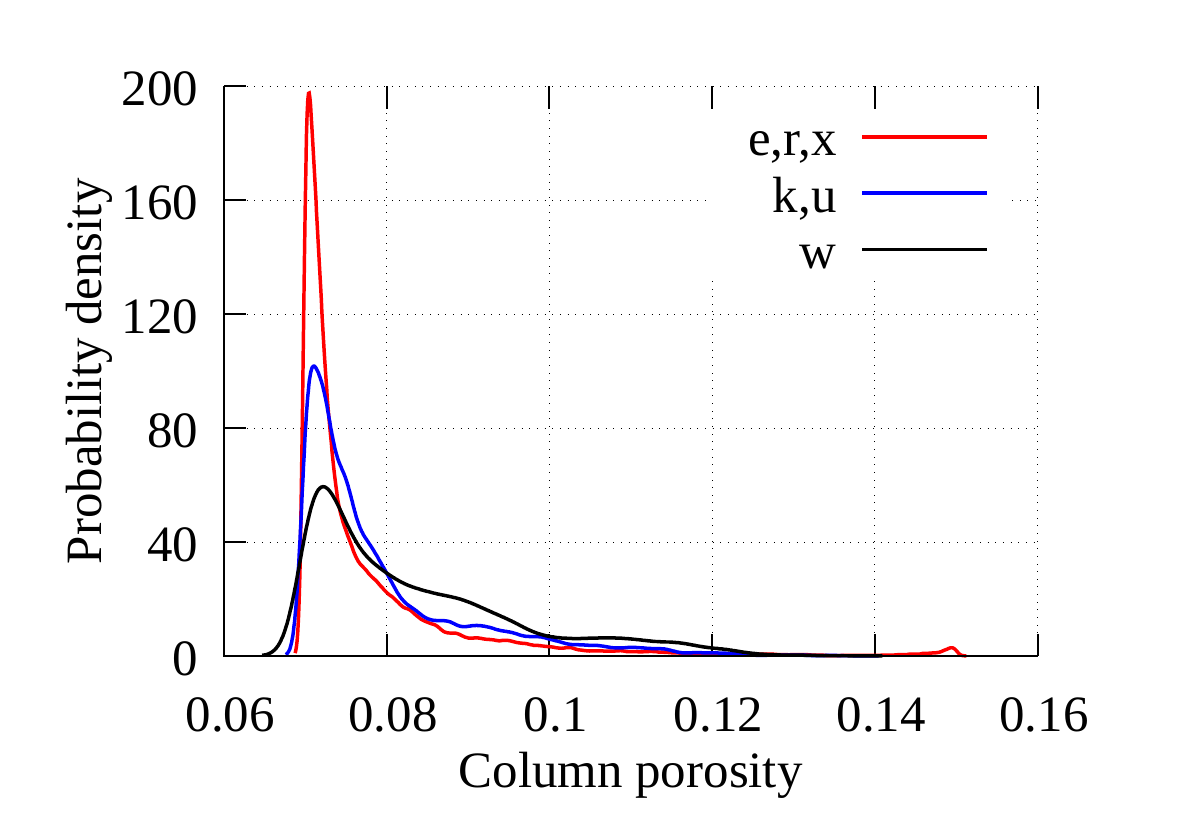}} \caption{Column porosity} \label{fig:postdist_col_b} \end{subfigure}
    \begin{subfigure}{0.33\textwidth}{\includegraphics[width=\textwidth]{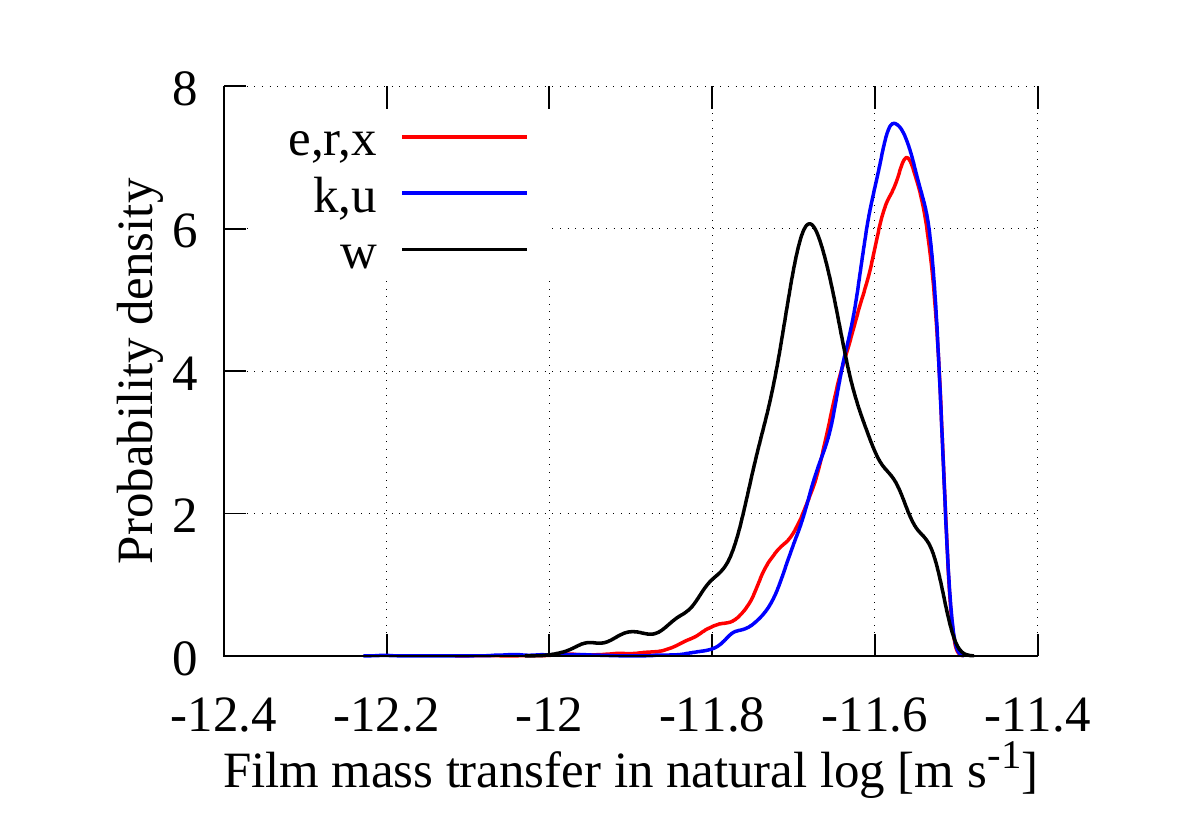}} \caption{Film mass transfer (log)} \label{fig:postdist_col_c} \end{subfigure}
    \begin{subfigure}{0.33\textwidth}{\includegraphics[width=\textwidth]{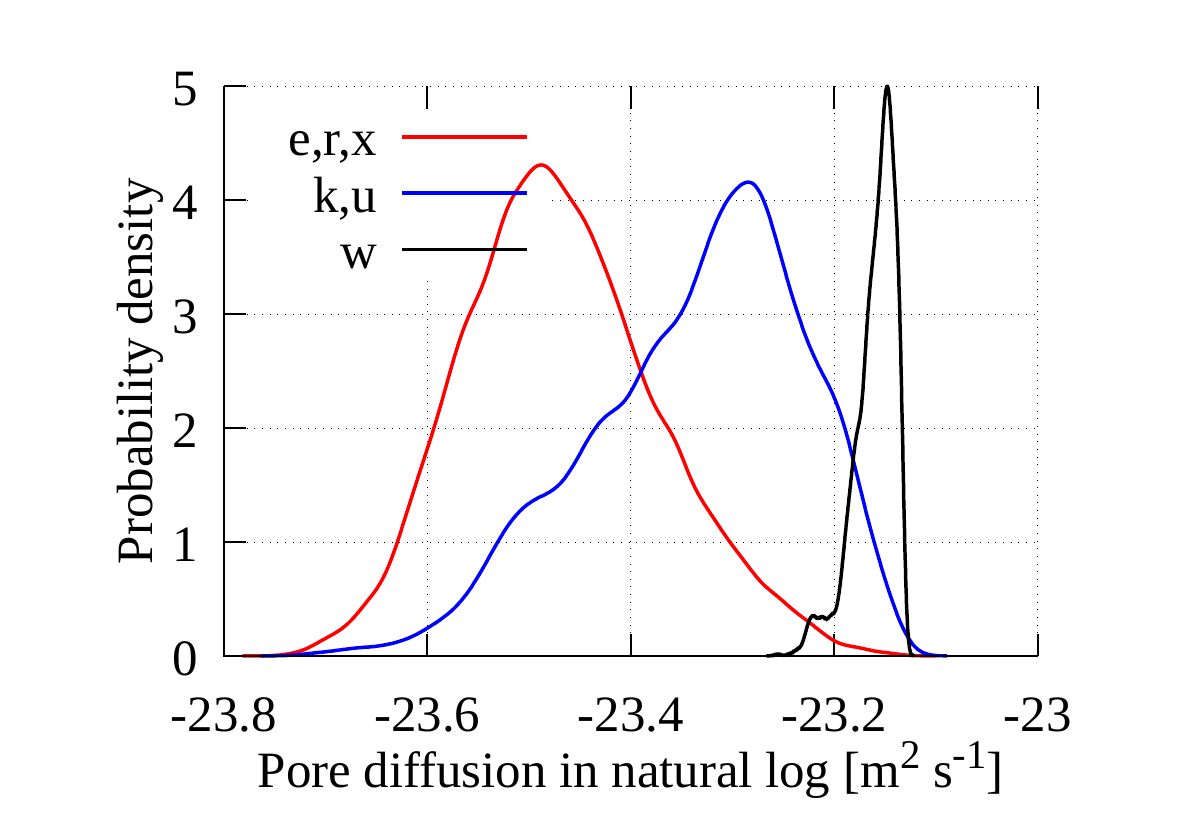}} \caption{Pore diffusion (log)} \label{fig:postdist_col_d} \end{subfigure}
    \begin{subfigure}{0.33\textwidth}{\includegraphics[width=\textwidth]{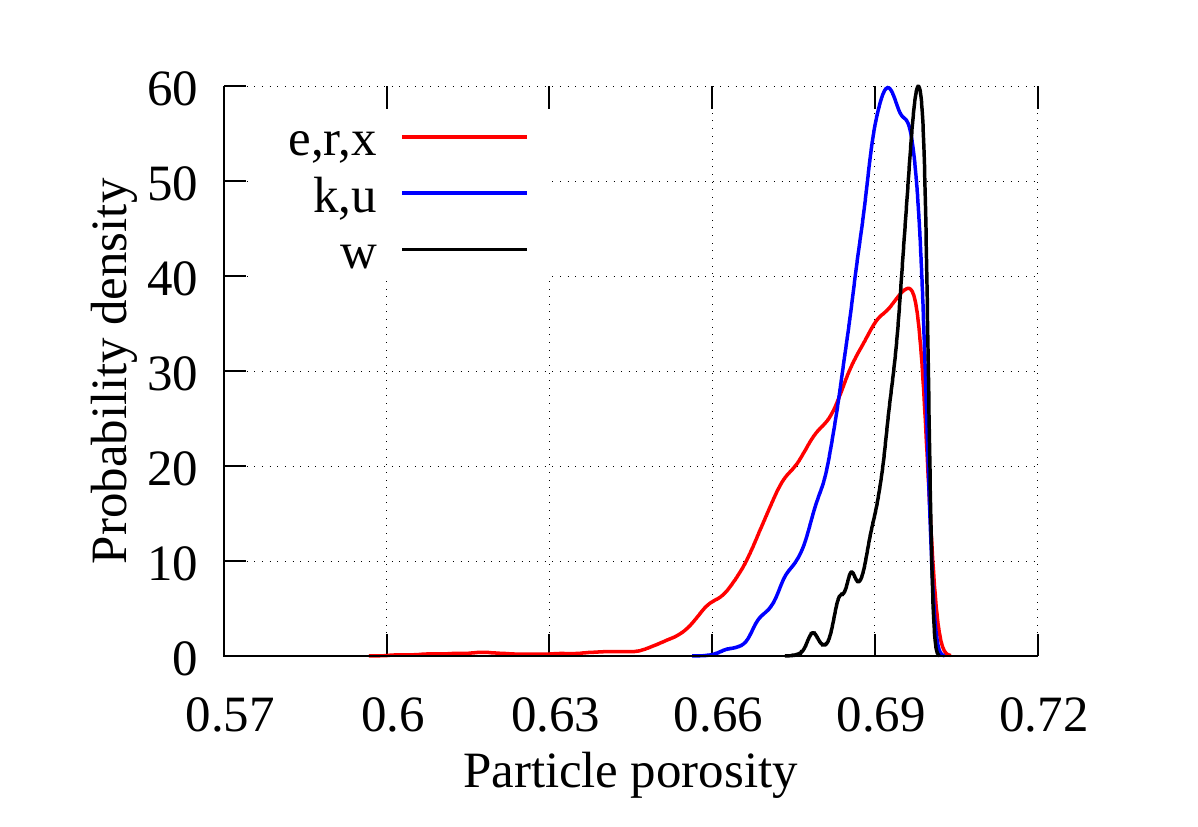}} \caption{Particle porosity} \label{fig:postdist_col_e} \end{subfigure}
    \begin{subfigure}{0.33\textwidth}{\includegraphics[width=\textwidth]{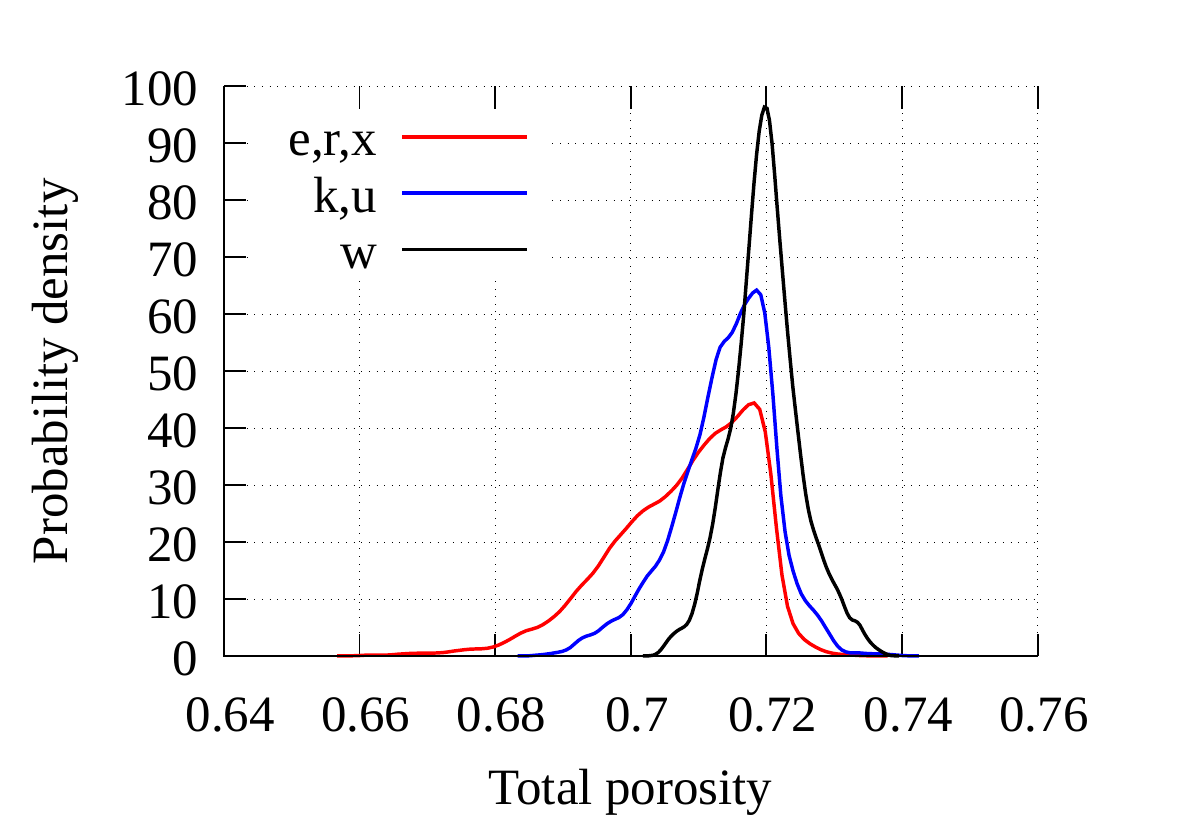}} \caption{Total porosity} \label{fig:postdist_col_f} \end{subfigure}
    \caption{Fourth stage results (separate analysis), binding properties from lysozyme step and gradient elution. Results of simultaneous analysis are shown for comparison (all). Posterior distributions in range plot mode (a-k), and total porosity (l).} 
    \label{fig:postdist_col}
\end{figure}

Column properties are characterized with model parameters of column porosity, $\varepsilon_c$, axial dispersion, $D_\text{ax}$, while particle properties by particle porosity, $\varepsilon_p$, pore diffusion, $D_p$ and film mass transfer, $k_f$.
For the column properties, distributions of $K=3$ modes still agree with each other.
However, for the particle properties, inconsistency among modes in the distribution of pore diffusion, $D_p$, is observed.
Thus, identification of particle properties could result in multi-modality.

Judging from the shapes of the distributions (cf.~Fig.~\ref{fig:postdist_col}), they are well-determined. 
Specifically, dispersion of the interstitial volume (ca. $D_\text{ax} =\SI{4.0e-10}{\square\metre\per\second}$) is larger than that of the porous volume (ca. $D_p =\SI{7.0e-11}{\square\metre\per\second}$), which meets our expectation.
The film mass transfer coefficient is peaked aroud \SI{1e-5}{\metre\per\second}; the particle porosity is peaked around \num{0.69}.
The coefficient values do generally agree with literature, even though data for the same stationary phase and flow conditions is unavailable.
For \textsc{sp} Sepharose \textsc{ff}  a value of $D_\text{ax} = \SI{1.57e-10}{\square\metre\per\second}$ has been published \citep{osberghaus2012determination}.
A pore diffusion coefficient, $D_p = \SI{6.07e-11}{\square\metre\per\second}$ and film mass transfer of $k_f = \SI{6.90e-6}{\metre\per\second}$ in \cite{puttmann2013fast}.
The distribution of column porosity, $\varepsilon_c$, peaked around \num{0.07} and with extremely low possibilities on \num{0.15}, is unexpected and abnormal, as it is way smaller than the normal values, \eg, \num{0.20}.
\cite{traub2005preparative} has experimentally reported a range between $[0.50, 0.90]$ for $\varepsilon_p$, a range between $[0.26, 0.48]$ for $\varepsilon_c$.
This can be partially attributed to the mirror implementation of \textsc{dpfr} and \textsc{cstr} after the column, which might occupy more porosities from the column side.
Running simulations without the identical $L^\textsc{dpfr}$ value can validate it.
Smaller values can also be potentially explained by mechanical compression during column packing and salt elution.
The value of total porosity, $\varepsilon_t$, is also slightly smaller than the reported values, see Fig.~\ref{fig:postdist_col_f}.

\subsubsection{Binding properties}
\begin{figure}[h]
    \centering
    \begin{subfigure}{0.32\textwidth}{\includegraphics[width=\textwidth]{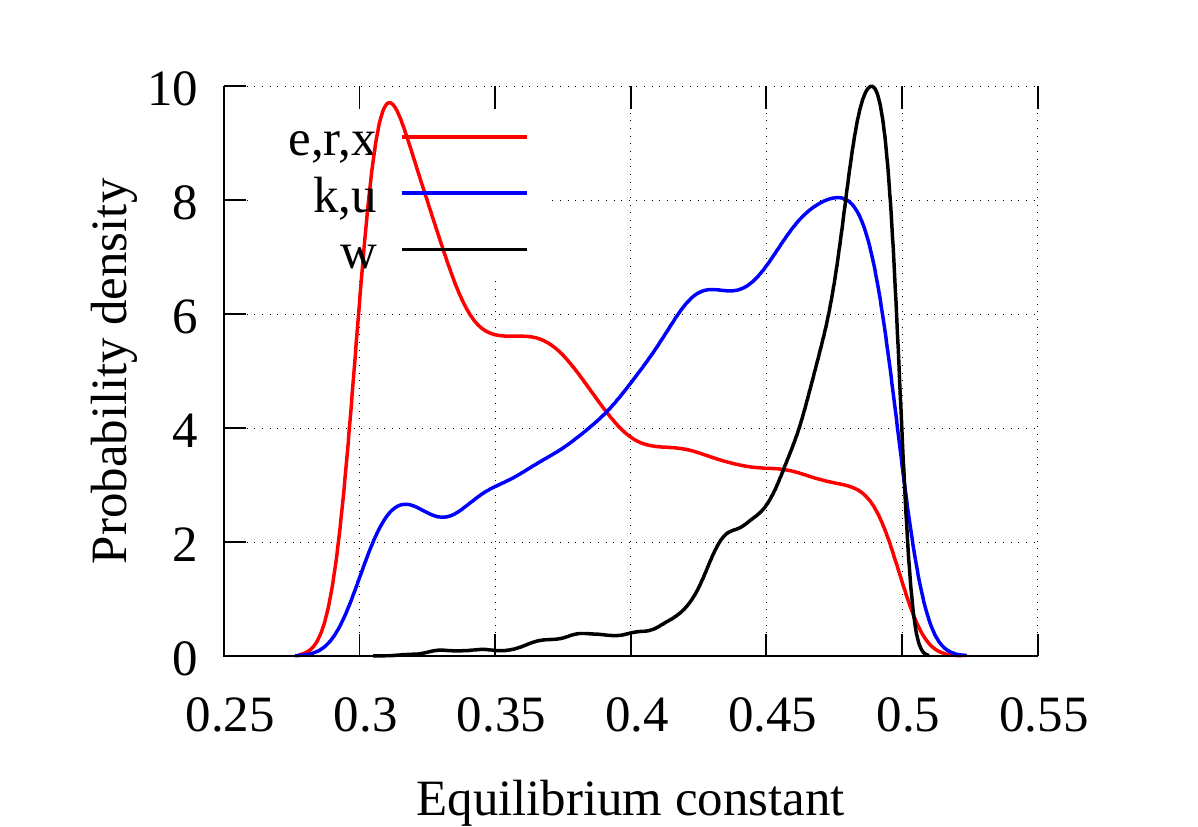}} \caption{Equilibrium constant} \label{fig:postdist_bind_a} \end{subfigure}
    \begin{subfigure}{0.32\textwidth}{\includegraphics[width=\textwidth]{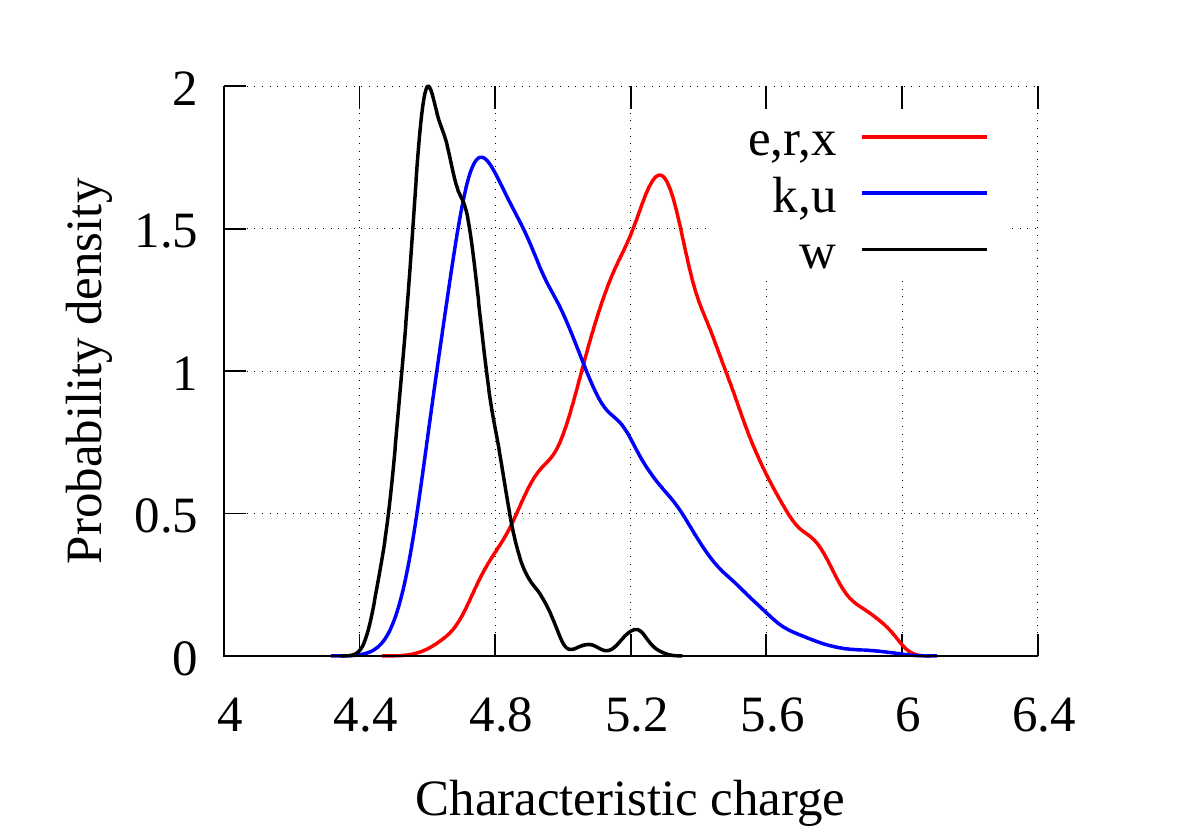}} \caption{Characteristic charge} \label{fig:postdist_bind_b} \end{subfigure}
    \begin{subfigure}{0.32\textwidth}{\includegraphics[width=\textwidth]{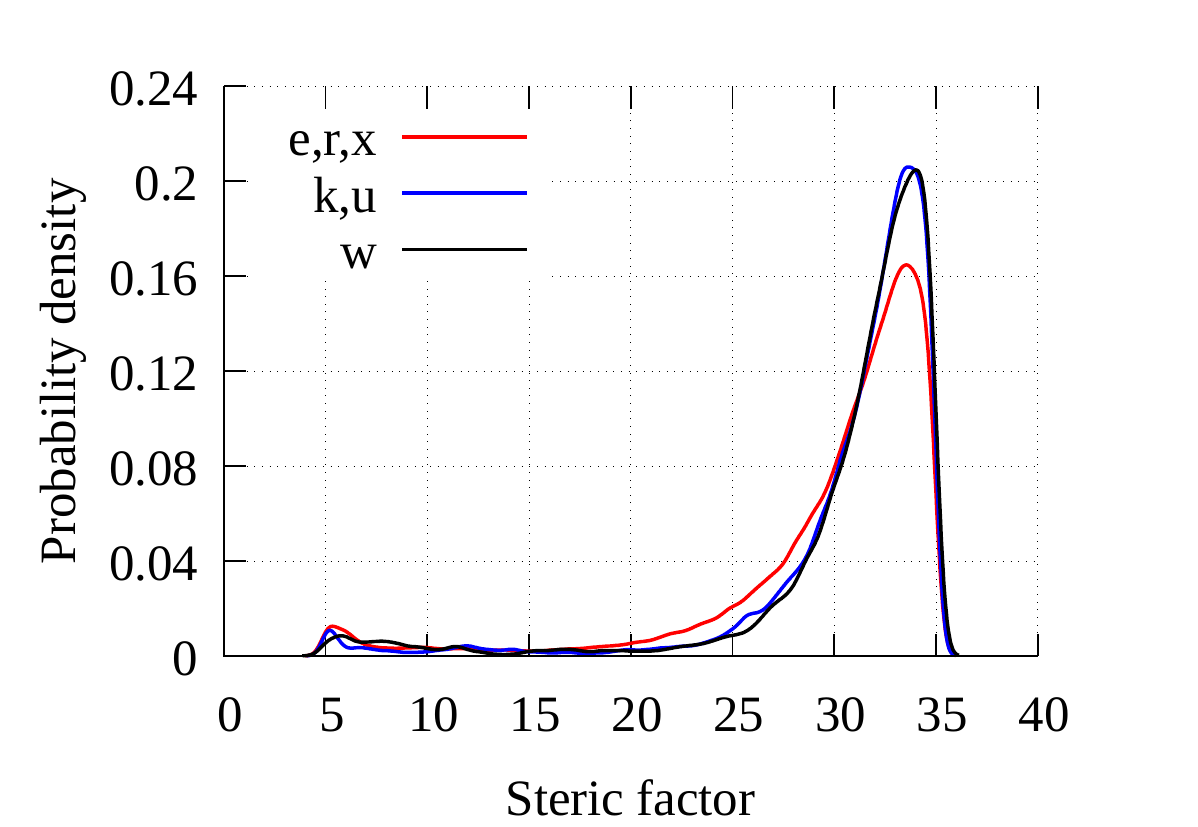}} \caption{Steric factor} \label{fig:postdist_bind_c} \end{subfigure}
    \caption{} 
    \label{fig:postdist_bind}
\end{figure}

Binding behaviour of the component onto the stationary phase is described with the \textsc{sma} isotherm.
The binding parameters at this stage, $k_\text{eq}$, $\nu$, $\sigma$, are not well-identified, as $K=3$ modes renders three different distributions of parameters $k_\text{eq}$ and $\nu$.
Notably, the distribution of steric factor has two modes, with one mode peaked around \num{5.0} and another mode around \num{35}.
\cite{osberghaus2012model} have reported values spread over a wide range from $\sigma= 20$ to $\sigma=50$, which agrees with our results.
The equilibrium constant, $k_\text{eq}$, has most probable value between \num{0.3} and \num{0.5}; this slightly differs with the literature values, $k_\text{eq} = 0.04$ \citep{osberghaus2012determination}, $k_\text{eq} = 0.14$ \citep{osberghaus2012optimizing} and $k_\text{eq} = 0.27$ \citep{osberghaus2012model}, reported for lysozyme on \textsc{sp} Sepharose \textsc{ff} at pH at \num{7.0} using the inverse methods.
Values of $\nu = 3.37$ \citep{osberghaus2012model}, $\nu = 3.40$ \citep{osberghaus2012optimizing} and $\nu = 4.72$ \citep{osberghaus2012determination} have been published for lysozyme on \textsc{sp} Sepharose \textsc{ff} at pH \num{7.0}, which agrees well with our results.

Therefore, binding parameters are relatively harder to identify in comparison with other model parameters, and have higher probability results in multi-modality in model calibration or parameter estimation.
The long tailing of steric factor (this observation follows our general experience that it is often poorly determined), $\sigma$, could be a clue of re-orientation of lysozyme.

\section{Conclusion}
Bayesian inference and \textsc{mcmc} algorithm are promising in rendering posterior distributions of model parameters upon uncertainties from experiments and models.
Difficulties in convergence diagnostic resulted from multi-modality of the target distribution can be solved with the proposed clustering-based $\widehat{R}$ diagnostic.

Both partitional ($K$-means and its variants, \eg, $K$-medians, $K$-medoids) and hierarchical (dendrogram representation) clustering methods have been tested and cross-validated with each other in this study, and consistent results have been observed.
The consistency of clustering among various methods might be attributed to the mapping implementation from model parameter values to smooth kernel density functions.
Such implementation extracts distinct features of the model parameter space, such as the Perron cluster analysis method that leverages eigen vales.

By utilizing the partitional clustering methods, it is challenging to estimate the correct number of clusters ($K$).
Researchers have proposed prominent methods for addressing the challenge; they would be considered in our later work.
It is straightforward and convenient to extend the clustering methods used in this study to other methods, \eg, the Perron cluster analysis.
Due to the computing capacity of our node, rather small parallel chains had been used.
It would be more powerful to utilize more parallel \textsc{mcmc} sampling chains or \textsc{gpu} source.
We would be rather positive on robustness of the clustering methods used in this study on large number of parallel chains.

\section{Acknowledgements}

\bibliographystyle{model5-names}
\bibliography{references.bib}
\appendix

\section{Finding optimal $K$ value}\label{sec:appd}
\begin{figure}[h]
    \centering
    \includegraphics[width=0.5\textwidth]{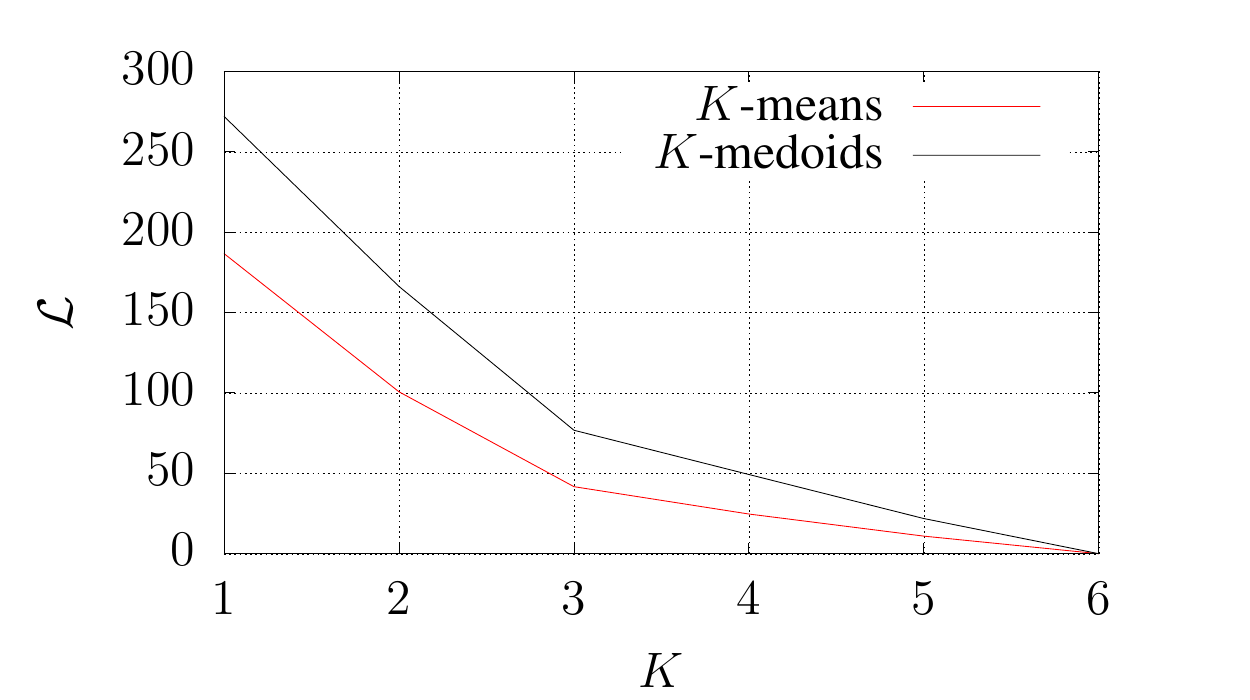}
    \caption{Distortion values over linear cluster values of two partitional clustering methods.}
    \label{fig:KL}
\end{figure}

In this study, the optimal $K$ value is determined firstly by running the $K$-means algorithm over a linear range of the $K$ values (\ie, from 1 to $p=6$).
Unarguably, this method is raw, but it works so far.
Thereafter, the elbow method is used to detect the optimal $K$ value from the $K$ vs.~$\mathcal{L}$ plot, see Fig.~\ref{fig:KL}.

The distortion values are calculated from Eq.~\eqref{eq:distortion}; although different clustering methods were used, the squared Euclidean distance was the same.
As shown in Fig.~\ref{fig:KL}, the optimal $K$ value can be claimed to be 3 from both the $K$-means and the $K$-medoids methods.
Subroutines of \texttt{kmeans} and \texttt{kmedoids} from the \textsc{matlab} were used.
The elbow searching method used in this study is taken from the MathWorks center, see \url{https://www.mathworks.com/matlabcentral/fileexchange/35094-knee-point} for more detailed information.

\end{document}